\documentclass[11pt]{article}
\usepackage{amsmath, upgreek}
\usepackage{amssymb}
\usepackage{color}
\usepackage{amscd}
\usepackage{xspace}
\usepackage{verbatim}
\usepackage{graphicx}
\renewcommand{\theequation}{\thesection.\arabic{equation}
}
\newcommand{\mysection}[1]{
\section{#1}\setcounter{equation}{0}}
\title{\bf Schr\"odinger operators   
with Leray-Hardy potential singular on the boundary}
\author{{\bf Huyuan Chen\footnote{\noindent Department of Mathematics, Jiangxi Normal University,
Nanchang 330022, China. E-mail: chenhuyuan@yeah.net}} \\[4mm]
 {\bf Laurent V\'eron \footnote{\noindent
Laboratoire de Math\'{e}matiques et Physique Th\'{e}orique, Universit\'e de Tours, 37200 Tours, France. E-mail: veronl@univ-tours.fr}}
}

\date{}
\begin{document}
 \maketitle


\newcommand{\txt}[1]{\;\text{ #1 }\;}
\newcommand{\tbf}{\textbf}
\newcommand{\tit}{\textit}
\newcommand{\tsc}{\textsc}
\newcommand{\trm}{\textrm}
\newcommand{\mbf}{\mathbf}
\newcommand{\mrm}{\mathrm}
\newcommand{\bsym}{\boldsymbol}
\newcommand{\scs}{\scriptstyle}
\newcommand{\sss}{\scriptscriptstyle}
\newcommand{\txts}{\textstyle}
\newcommand{\dsps}{\displaystyle}
\newcommand{\fnz}{\footnotesize}
\newcommand{\scz}{\scriptsize}
\newcommand{\be}{\begin{equation}}
\newcommand{\bel}[1]{\begin{equation}\label{#1}}
\newcommand{\ee}{\end{equation}}
\newcommand{\eqnl}[2]{\begin{equation}\label{#1}{#2}\end{equation}}
\newcommand{\barr}{\begin{eqnarray}}
\newcommand{\earr}{\end{eqnarray}}
\newcommand{\bars}{\begin{eqnarray*}}
\newcommand{\ears}{\end{eqnarray*}}
\newcommand{\nnu}{\nonumber \\}
\newtheorem{subn}{\name}
\renewcommand{\thesubn}{}
\newcommand{\bsn}[1]{\def\name{#1}\begin{subn}}
\newcommand{\esn}{\end{subn}}
\newtheorem{sub}{\name}[section]
\newcommand{\dn}[1]{\def\name{#1}}   
\newcommand{\bs}{\begin{sub}}
\newcommand{\es}{\end{sub}}
\newcommand{\bsl}[1]{\begin{sub}\label{#1}}
\newcommand{\bth}[1]{\def\name{Theorem}
\begin{sub}\label{t:#1}}
\newcommand{\blemma}[1]{\def\name{Lemma}
\begin{sub}\label{l:#1}}
\newcommand{\bcor}[1]{\def\name{Corollary}
\begin{sub}\label{c:#1}}
\newcommand{\bdef}[1]{\def\name{Definition}
\begin{sub}\label{d:#1}}
\newcommand{\bprop}[1]{\def\name{Proposition}
\begin{sub}\label{p:#1}}

\newcommand{\aand}{\quad\mbox{and}\quad}
\newcommand{\M}{{\cal M}}
\newcommand{\A}{{\cal A}}
\newcommand{\B}{{\cal B}}
\newcommand{\I}{{\cal I}}
\newcommand{\J}{{\cal J}}
\newcommand{\D}{\displaystyle}
\newcommand{\RR}{ I\!\!R}
\newcommand{\C}{\mathbb{C}}
\newcommand{\R}{\mathbb{R}}
\newcommand{\Z}{\mathbb{Z}}
\newcommand{\N}{\mathbb{N}}
\newcommand{\T}{{\rm T}^n}
\newcommand{\cuad}{{\sqcap\kern-.68em\sqcup}}
\newcommand{\abs}[1]{\mid #1 \mid}
\newcommand{\norm}[1]{\|#1\|}
\newcommand{\equ}[1]{(\ref{#1})}
\newcommand\rn{\mathbb{R}^N}
\renewcommand{\theequation}{\thesection.\arabic{equation}}
\newtheorem{definition}{Definition}[section]
\newtheorem{theorem}{Theorem}[section]
\newtheorem{proposition}{Proposition}[section]
\newtheorem{example}{Example}[section]
\newtheorem{proof}{proof}[section]
\newtheorem{lemma}{Lemma}[section]
\newtheorem{corollary}{Corollary}[section]
\newtheorem{remark}{Remark}[section]
\newcommand{\bremark}{\begin{remark} \em}
\newcommand{\eremark}{\end{remark} }
\newtheorem{claim}{Claim}


\newcommand{\rth}[1]{Theorem~\ref{t:#1}}
\newcommand{\rlemma}[1]{Lemma~\ref{l:#1}}
\newcommand{\rcor}[1]{Corollary~\ref{c:#1}}
\newcommand{\rdef}[1]{Definition~\ref{d:#1}}
\newcommand{\rprop}[1]{Proposition~\ref{p:#1}}
\newcommand{\BA}{\begin{array}}
\newcommand{\EA}{\end{array}}
\newcommand{\BAN}{\renewcommand{\arraystretch}{1.2}
\setlength{\arraycolsep}{2pt}\begin{array}}
\newcommand{\BAV}[2]{\renewcommand{\arraystretch}{#1}
\setlength{\arraycolsep}{#2}\begin{array}}
\newcommand{\BSA}{\begin{subarray}}
\newcommand{\ESA}{\end{subarray}}
\newcommand{\BAL}{\begin{aligned}}
\newcommand{\EAL}{\end{aligned}}
\newcommand{\BALG}{\begin{alignat}}
\newcommand{\EALG}{\end{alignat}}
\newcommand{\BALGN}{\begin{alignat*}}
\newcommand{\EALGN}{\end{alignat*}}
\newcommand{\note}[1]{\textit{#1.}\hspace{2mm}}
\newcommand{\Proof}{\note{Proof}}
\newcommand{\qeda}{\hspace{10mm}\hfill $\square$}
\newcommand{\qed}{\\
${}$ \hfill $\square$}
\newcommand{\Remark}{\note{Remark}}
\newcommand{\modin}{$\,$\\[-4mm] \indent}
\newcommand{\forevery}{\quad \forall}
\newcommand{\set}[1]{\{#1\}}
\newcommand{\setdef}[2]{\{\,#1:\,#2\,\}}
\newcommand{\setm}[2]{\{\,#1\mid #2\,\}}
\newcommand{\mt}{\mapsto}
\newcommand{\lra}{\longrightarrow}
\newcommand{\lla}{\longleftarrow}
\newcommand{\llra}{\longleftrightarrow}
\newcommand{\Lra}{\Longrightarrow}
\newcommand{\Lla}{\Longleftarrow}
\newcommand{\Llra}{\Longleftrightarrow}
\newcommand{\warrow}{\rightharpoonup}
\newcommand{
\paran}[1]{\left (#1 \right )}
\newcommand{\sqbr}[1]{\left [#1 \right ]}
\newcommand{\curlybr}[1]{\left \{#1 \right \}}
\newcommand{
\paranb}[1]{\big (#1 \big )}
\newcommand{\lsqbrb}[1]{\big [#1 \big ]}
\newcommand{\lcurlybrb}[1]{\big \{#1 \big \}}
\newcommand{\absb}[1]{\big |#1\big |}
\newcommand{\normb}[1]{\big \|#1\big \|}
\newcommand{
\paranB}[1]{\Big (#1 \Big )}
\newcommand{\absB}[1]{\Big |#1\Big |}
\newcommand{\normB}[1]{\Big \|#1\Big \|}
\newcommand{\produal}[1]{\langle #1 \rangle}

\newcommand{\thkl}{\rule[-.5mm]{.3mm}{3mm}}
\newcommand{\thknorm}[1]{\thkl #1 \thkl\,}
\newcommand{\trinorm}[1]{|\!|\!| #1 |\!|\!|\,}
\newcommand{\bang}[1]{\langle #1 \rangle}
\def\angb<#1>{\langle #1 \rangle}
\newcommand{\vstrut}[1]{\rule{0mm}{#1}}
\newcommand{\rec}[1]{\frac{1}{#1}}
\newcommand{\opname}[1]{\mbox{\rm #1}\,}
\newcommand{\supp}{\opname{supp}}
\newcommand{\dist}{\opname{dist}}
\newcommand{\myfrac}[2]{{\displaystyle \frac{#1}{#2} }}
\newcommand{\myint}[2]{{\displaystyle \int_{#1}^{#2}}}
\newcommand{\mysum}[2]{{\displaystyle \sum_{#1}^{#2}}}
\newcommand {\dint}{{\displaystyle \myint\!\!\myint}}
\newcommand{\q}{\quad}
\newcommand{\qq}{\qquad}
\newcommand{\hsp}[1]{\hspace{#1mm}}
\newcommand{\vsp}[1]{\vspace{#1mm}}
\newcommand{\ity}{\infty}
\newcommand{\prt}{\partial}
\newcommand{\sms}{\setminus}
\newcommand{\ems}{\emptyset}
\newcommand{\ti}{\times}
\newcommand{\pr}{^\prime}
\newcommand{\ppr}{^{\prime\prime}}
\newcommand{\tl}{\tilde}
\newcommand{\sbs}{\subset}
\newcommand{\sbeq}{\subseteq}
\newcommand{\nind}{\noindent}
\newcommand{\ind}{\indent}
\newcommand{\ovl}{\overline}
\newcommand{\unl}{\underline}
\newcommand{\nin}{\not\in}
\newcommand{\pfrac}[2]{\genfrac{(}{)}{}{}{#1}{#2}}

\def\ga{\alpha}     \def\gb{\beta}       \def\gg{\gamma}
\def\gc{\chi}       \def\gd{\delta}      \def\ge{\epsilon}
\def\gth{\theta}                         \def\vge{\varepsilon}
\def\gf{\phi}       \def\vgf{\varphi}    \def\gh{\eta}
\def\gi{\iota}      \def\gk{\kappa}      \def\gl{\lambda}
\def\gm{\mu}        \def\gn{\nu}         \def\gp{\pi}
\def\vgp{\varpi}    \def\gr{\rho}        \def\vgr{\varrho}
\def\gs{\sigma}     \def\vgs{\varsigma}  \def\gt{\tau}
\def\gu{\upsilon}   \def\gv{\vartheta}   \def\gw{\omega}
\def\gx{\xi}        \def\gy{\psi}        \def\gz{\zeta}
\def\Gg{\Gamma}     \def\Gd{\Delta}      \def\Gf{\Phi}
\def\Gth{\Theta}
\def\Gl{\Lambda}    \def\Gs{\Sigma}      \def\Gp{\Pi}
\def\Gw{\Omega}     \def\Gx{\Xi}         \def\Gy{\Psi}

\def\CS{{\mathcal S}}   \def\CM{{\mathcal M}}   \def\CN{{\mathcal N}}
\def\CR{{\mathcal R}}   \def\CO{{\mathcal O}}   \def\CP{{\mathcal P}}
\def\CA{{\mathcal A}}   \def\CB{{\mathcal B}}   \def\CC{{\mathcal C}}
\def\CD{{\mathcal D}}   \def\CE{{\mathcal E}}   \def\CF{{\mathcal F}}
\def\CG{{\mathcal G}}   \def\CH{{\mathcal H}}   \def\CI{{\mathcal I}}
\def\CJ{{\mathcal J}}   \def\CK{{\mathcal K}}   \def\CL{{\mathcal L}}
\def\CT{{\mathcal T}}   \def\CU{{\mathcal U}}   \def\CV{{\mathcal V}}
\def\CZ{{\mathcal Z}}   \def\CX{{\mathcal X}}   \def\CY{{\mathcal Y}}
\def\CW{{\mathcal W}} \def\CQ{{\mathcal Q}}
\def\BBA {\mathbb A}   \def\BBb {\mathbb B}    \def\BBC {\mathbb C}
\def\BBD {\mathbb D}   \def\BBE {\mathbb E}    \def\BBF {\mathbb F}
\def\BBG {\mathbb G}   \def\BBH {\mathbb H}    \def\BBI {\mathbb I}
\def\BBJ {\mathbb J}   \def\BBK {\mathbb K}    \def\BBL {\mathbb L}
\def\BBM {\mathbb M}   \def\BBN {\mathbb N}    \def\BBO {\mathbb O}
\def\BBP {\mathbb P}   \def\BBR {\mathbb R}    \def\BBS {\mathbb S}
\def\BBT {\mathbb T}   \def\BBU {\mathbb U}    \def\BBV {\mathbb V}
\def\BBW {\mathbb W}   \def\BBX {\mathbb X}    \def\BBY {\mathbb Y}
\def\BBZ {\mathbb Z}

\def\GTA {\mathfrak A}   \def\GTB {\mathfrak B}    \def\GTC {\mathfrak C}
\def\GTD {\mathfrak D}   \def\GTE {\mathfrak E}    \def\GTF {\mathfrak F}
\def\GTG {\mathfrak G}   \def\GTH {\mathfrak H}    \def\GTI {\mathfrak I}
\def\GTJ {\mathfrak J}   \def\GTK {\mathfrak K}    \def\GTL {\mathfrak L}
\def\GTM {\mathfrak M}   \def\GTN {\mathfrak N}    \def\GTO {\mathfrak O}
\def\GTP {\mathfrak P}   \def\GTR {\mathfrak R}    \def\GTS {\mathfrak S}
\def\GTT {\mathfrak T}   \def\GTU {\mathfrak U}    \def\GTV {\mathfrak V}
\def\GTW {\mathfrak W}   \def\GTX {\mathfrak X}    \def\GTY {\mathfrak Y}
\def\GTZ {\mathfrak Z}   \def\GTQ {\mathfrak Q}

\font\Sym= msam10 
\def\SYM#1{\hbox{\Sym #1}}
\newcommand{\bdw}{\prt\Gw\xspace}
\date{}
\maketitle

\begin{abstract}
We study the kernel function of the operator $u\mapsto \CL_\gm u=-\Gd u+\tfrac {\gm}{|x|^2}u$ in a bounded smooth domain $\Gw\subset\BBR^N_+$ such that $0\in\prt\Gw$, where $\gm\geq-\frac{N^2}{4}$ is a constant. We show the existence of a Poisson kernel vanishing at $0$ and a singular kernel with a singularity at $0$.  We prove the existence and uniqueness of weak solutions of $\CL_\gm u=0$ in $\Gw$ with boundary data $\gn+k\gd_0$, where $\gn$ is a Radon measure on $\prt\Gw\setminus\{0\}$, $k\in\R$ and show that this boundary data corresponds in a unique way to the boundary trace of positive solution of $\CL_\gm u=0$ in $\Gw$.
\end{abstract}

\noindent
  \noindent {\small {\bf Key Words}:   Hardy Potential, Harnack inequality, limit set, Radon Measure.  }\vspace{1mm}

\noindent {\small {\bf MSC2010}:  35B44, 35J75. }\tableofcontents
\vspace{1mm}
\hspace{.05in}
\medskip

\setcounter{equation}{0}
\section{Introduction}
We denote by $\CL_\gm$ the  Schr\"odinger operator defined in a domain $\Gw\subset\BBR^N$ by
$$
\CL_\gm u:=-\Gd u+\frac{\gm}{|x|^2}u,
$$
where $\gm$ is a real constant and $N\geq 2$. This operator which is associated to the Hardy inequality has been thoroughly studied in the last thirty years. When the singular point $0$ belongs to $\Gw$, it appears a critical value 
 $$
\gm_0=-\left(\myfrac{N-2}{2}\right)^2
$$
and the range of the $\gm$ in which the operator is bounded from below is $[\gm_0,\infty)$. This is linked to the Hardy inequality
$$
\myint{\Gw}{}|\nabla \phi|^2+\gm_0\myint{\Gw}{}\myfrac{\phi^2}{|x|^2}dx\geq 0\quad\text{for all }\gf\in C^\infty_0(\Gw).
$$
Furthermore this inequality is never achieved if $\Gw$ is bounded, in which case a remainder was shown to exist by Br\'ezis and V\'azquez  \cite{BrVa}. When $\gl$ is a Radon measure in $\Gw$, the associated Dirichlet problem 
$$
\left\{\BA{lll}
\CL_\gm u=\gl&\quad\text{in }\ \Gw,\\[0.5mm]
\phantom{\CL_\gm}
u=0&\quad\text{on } \;\prt\Gw
\EA\right.$$
is studied in its full generality in \cite{CQZ} and \cite{ChVe} thanks to the introduction of a notion of very weak solution associated to 
some specific weight. Thanks to this new formulation an extensive treatment of the associated semilinear problem
$$
 \left\{\BA {lll}
\CL_\gm u+ g(u)=\gl&\quad\text{in }\ \Gw,\\[0.5mm]
\phantom{\CL_\gm + g(u)}
u=0&\quad\text{on }\, \prt\Gw,
\EA\right.$$
where $g:\BBR\mapsto\BBR$ is a continuous nondecreasing function is developed in \cite{ChVe}. \smallskip

In this article we assume that the singular point of the potential lies on the boundary of the domain  $\Gw$, and we are mainly interested in the two problems:\\
 1- To define a notion of  {\it very weak solution} for the problem
 \bel{1-A5} \left\{ \BA {lll}
\CL_\gm u=0&\quad\text{in }\, \Gw,\\[0.5mm]
\phantom{\CL_\gm}
u=\gn&\quad\text{on }\, \prt\Gw,
\EA\right.\ee
where $\gn$ is a Radon measure on $\prt\Gw$, and more generaly on $\prt\Gw\setminus\{0\}$;\\
 2- To prove the existence of a boundary trace for any positive {\it $\CL_\gm$-harmonic function}, i.e. solution of $\CL_\gm u=0$ in $\Gw$
 and to connect it to the problem (\ref{1-A5}).\smallskip

The model example is 
$\Gw=\BBR^N_+:=\{x=(x',x_N)\in\R^{N-1}\times\R:\, x_N>0\}$ although it is not a bounded domain. There exists a critical value 
\bel{1-A11}
\gm\geq\gm_1:=-\frac {N^2}{4}. 
\ee
This value is fundamental for the operator $\CL_\gm$ to be bounded from below since there holds,
\bel{1-Ax12}
\myint{\BBR^{N}_+}{}|\nabla \phi|^2+\gm_1\myint{\BBR^{N}_+}{}\myfrac{\phi^2}{|x|^2}dx\geq 0\quad\text{for all }\gf\in C^\infty_0(\BBR^{N}_+).
\ee 
The analysis of the model case is explicit. Let $(r,\gs)\in \BBR_+\ti \BBS^{N-1}_+$ be the spherical coordinates in $\BBR^N_+$, then, if (\ref{1-A11}) is satisfied, there exist two different types of positive $\CL_\gm$-harmonic functions vanishing on $\prt\BBR^{N}_+\setminus\{0\}$),
\bel{1-A13}
\gg_\gm(r,\gs)=r^{\ga_+}\psi_1(\gs)\quad\text{and }\;\gf_\gm(r,\gs)=\left\{\BA {lll}r^{\ga_-}\psi_1(\gs)\quad&\text{if }\gm>\gm_1,\\[2mm]
r^{-\frac{N-2}{2}}\ln(r^{-1})\psi_1(\gs)\quad&\text{if }\gm=\gm_1,\EA\right.
\ee 
where $\psi_1(\gs)=\frac{x_N}{|x|}$ generates $\ker (-\Gd'+(N-1)I)$ in $H^1_0(\BBS^{N-1}_+)$, and where
\bel{1-A14}
\ga_+:=\ga_+(\gm)=\frac{2-N}{2}+\sqrt{\gm+\frac{N^2}{4}}\quad\text{and }\;\ga_-:=\ga_-(\gm)=\frac{2-N}{2}-\sqrt{\gm+\frac{N^2}{4}}.
\ee
Put $d\gg_\gm(x)=\gg_\gm(x) dx$. We define the $\gg_\gm$-dual operator $\CL^*_\gm$ of $\CL_\gm$ by
\bel{1-A17}
\CL^*_\gm\gz=-\Gd\gz-\myfrac{2}{\gg_{\gm}}\langle\nabla\gg_{\gm},\nabla\gz\rangle\quad\text{for all }\,\gz\in C^2(\overline\BBR^N_+),
\ee
and we prove that $\gf_\gm$ is, in some sense, the fundamental solution of 
 $$
 \left\{\BA {lll}
\CL_\gm u=0&\quad\text{in }\ \BBR^N_+,\\[0.5mm]
\phantom{\CL_\gm}
u=\gd_0&\quad\text{on }\,\prt\BBR^N_+,
\EA\right.$$
since it satisfies 
$$
\myint{\BBR^N_+}{}\gf_\gm\CL^*_\gm\gz d\gg_\gm(x)=c_\gm\gz(0)\quad \text{for all $\gz\in C_c(\overline{\BBR^N_+})\cap C^{1,1}({\BBR^N_+})$}
$$
   such that $\gr\CL^*_\gm\gz\in L^{\infty}(\BBR^N_+)$, where  $c_\gm>0$ is a normalized constant
and $\gr(x)=\text{dist}(x,\partial \Omega)$. Here $\rho(x)=x_N$ when $\Omega=\R^N_+$. \medskip

When $\BBR^{N}_+$ is replaced by a bounded domain $\Gw$ satisfying the condition 
$$(\CC\text{-}1) \qquad\qquad\qquad 0\in \prt\Gw\text{ , }\;\Gw\subset\BBR^{N}_+\text{ and }\;\langle x,{\bf n}\rangle=O(|x|^2)\,\text{ for all }\;x\in\prt \Gw,\qquad\qquad\qquad\qquad\qquad\qquad$$
where ${\bf n}={\bf n}_x$ is the outward normal vector at $x$, inequality (\ref{1-Ax12}) holds but it is never achieved in the Hilbert space $H^1_0(\Omega)$.  Note that the last condition in $(\CC\text{-}1)$ holds if $\Gw$ is a $C^2$ domain. It is proved in \cite{Caz1}  that there exists a remainder under the following form:  
\bel{1-A12}
\myint{\Gw}{}|\nabla \phi|^2+\gm_1\myint{\Gw}{}\myfrac{\phi^2}{|x|^2}dx\geq \myfrac{1}{4}
\myint{\Gw}{}\myfrac{\phi^2}{|x|^2\ln^2(|x|R_\Gw^{-1})}dx
\qquad\text{for all }\gf\in C^\infty_c(\Gw),
\ee
where
$\displaystyle R_\Gw=\max_{z\in\Gw} |z|$. 
 Under the assumption  $(\CC\text{-}1)$,   there holds 
$$
\ell_\gm^\Gw:=\inf\left\{\myint{\Gw}{}\left(|\nabla v|^2+\myfrac{\gm}{|x|^2}v^2\right)dx:v\in C^{1}_c(\Gw),\myint{\Gw}{}v^2 dx=1
\right\}>0.$$
This first eigenvalue is achieved in $H^1_0(\Gw)$ if $\gm>\gm_1$, or in the space $H(\Gw)$ which is the closure of $C^{1}_c(\Gw)$ for the norm
 $$
 v\mapsto \norm v_{H(\Gw)}:=\sqrt{\myint{\Gw}{}\left(|\nabla v|^2+\myfrac{\gm_1}{|x|^2}v^2\right)dx},
$$
when $\gm=\gm_1$. In the sequel we set 
 $$
 H_\gm(\Gw)=\left\{\BA {lll} H^1_0(\Gw)&\qquad\text{if }\,\gm>\gm_1\\[1mm]
 H(\Gw)&\qquad\text{if }\,\gm=\gm_1.
 \EA\right.
 $$
Moreover, under the assumption $(\CC\text{-}1)$ the imbedding of  $H_\gm(\Gw)$ in $L^2(\Gw)$ is compact (see e.g. \cite{Caz1-2}). We denote by  $\gg_\gm^\Gw$ the positive eigenfunction, its satisfies
\bel{1-A22}
 \left\{\BA {lll}
 \CL_\gm\gg_\gm^\Gw=\ell_\gm^\Gw\gg_\gm^\Gw\quad&\text{in }\ \Gw,\\[2mm]
 \phantom{\CL_\gm}
 \gg_\gm^\Gw=0\quad&\text{on }\, \prt\Gw\setminus\{0\}.
\EA\right. \ee
We prove that there exist $c_j=c_j(\Gw,\gm)>0$, j=1, 2,  such that 
  \bel{1-C1}\BA {lll}
(i)\qquad\qquad\qquad \gg_\gm^\Gw(x)= c_1\gr(x)|x|^{\ga_+-1}(1+o (1))\quad\text{as }\;x\to 0,\\[2mm]

(ii)\qquad\qquad\qquad  |\nabla \gg_\gm^\Gw(x)|\leq c_2\myfrac{\gg_\gm^\Gw(x)}{\gr(x)}\quad \text{for all }\; x\in\Gw.
\EA\ee
 This function will play the role of a weight function for replacing $\gg_\gm$.  Next we construct the Poisson kernel $\BBK^\Gw_\gm$ of $\CL_\gm$ in $\Gw\ti\prt\Gw$. When $\gm\geq 0$ this construction can be made by truncation as in \cite{VeYa}, considering for $\ge>0$ and $\gl\in\mathfrak M_+(\prt\Gw)$ the solution $u_\ge$ of 
$$
 \left\{\BA {lll}
-\Gd u+\myfrac{\gm}{\max\{\ge^2,|x|^2\}}u=0\quad&\text{in }\ \Gw,\\[1.5mm]
 \phantom{-\Gd u+\myfrac{\gm}{\max\{\ge^2,|x|^2\}}}
u=\gl\quad&\text{on }\, \prt\Gw.
\EA\right.$$
By a more elaborate method, we also construct the Poisson kernel when $\gm_1\leq\gm<0$. It is important to notice that when $\gm>0$ the kernel has the property that 
 \bel{1-A24}K^\Gw_\gm(x,0)=0\quad\text{for all }\,x\in\overline\Gw\setminus\{0\}
 \ee
by \cite[Theorem A.1]{VeYa}.  Because of (\ref{1-A24}) it is clear that the Poisson kernel cannot be the tool for describing all the positive $\CL_\gm$-harmonic functions.  
Our first concern in this article is to clarify the Poisson kernel of  $\CL_\gm$.

We first characterize the positive $\CL_\gm$-harmonic functions which are singular at $0$.\medskip

\nind{\bf Theorem A} {\it Let $\Gw$ be a $C^2$ bounded domain such that $0\in\prt\Gw$ and $\gm\geq\gm_1$. If $u$ is a nonnegative 
$\CL_\gm$-harmonic function vanishing on $B_{r_0}(0)\cap(\prt\Gw\setminus\{0\})$ for some $r_0>0$, there  exists $k\geq 0$ such that 
$$
\lim_{x\to 0}\myfrac{u(x)}{\gr(x)|x|^{\ga_--1}}=k,
$$
if $\gm>\gm_1$ and 
$$
\lim_{x\to 0}\myfrac{|x|^{\frac{N}{2}}u(x)}{\gr(x)\ln |x|}=-k,
$$
if $\gm=\gm_1$.
}
\medskip

Actually the above convergences hold in a stronger way. In order to prove that such solutions truly exist we construct the {\it kernel function} $\phi^\Gw_\gm$ (see  \cite{HuWh} for the denomination)  which is the analogue in a bounded domain of the explicit singular solution $\phi_\gm$ defined in $\BBR^N_+$. \medskip

\nind{\bf Theorem B} {\it Let $\Gw$ be a $C^2$ bounded domain such that $0\in\prt\Gw$ satisfying $(\CC\text{-}1)$ and $\gm\geq\gm_1$. Then there  exists a positive $\CL_\gm$-harmonic function in $\Gw$, which vanishes on $\partial\Gw\setminus\{0\}$ which satisfies,
  \bel{1-A25}\BA {lll}\displaystyle
\phi^\Gw_\gm(x)= \gr(x)|x|^{\ga_--1}(1+o (1))\quad\text{as }\;x\to 0,
\EA\ee
 if $\gm>\gm_1$, and 
  \bel{1-A26}\BA {lll}
\phi^\Gw_{\gm_1}(x)=\gr(x)|x|^{-\frac{N}{2}}(|\ln|x||+1)(1+o (1))\quad\text{as }\;x\to 0,
\EA\ee
 if $\gm=\gm_1$.  
}
\medskip

As in the model case, we define the $\gg_\gm^\Gw$-dual operator of $\CL_\gm$ by
 $$
\CL^*_\gm\gz=-\Gd\gz-\myfrac{2}{\gg_\gm^\Gw}\langle\nabla\gg_\gm^\Gw,\nabla\gz\rangle+\ell_\gm^\Gw\gz\quad\text{for all }\,\gz\in C^{1,1}(\Gw).
 $$
The following commutation formula holds
   \bel{2-0B0}\BA {lll}
\CL_\gm(\gg^\Gw_\gm\gz)=\gg^\Gw_\gm\CL^*_\gm\gz.
\EA\ee

\nind{\bf Corollary C} {\it Let $\Gw$ be a $C^2$ bounded domain such that $0\in\prt\Gw$ satisfying $(\CC$-$1)$ and $\gm\geq\gm_1$. Then $\gf_\gm^\Gw$ is the unique function belonging to $L^1(\Gw,\gr^{-1}d\gg_\gm^\Gw)$ which satisfies 
  \bel{1-Ax32z}
\myint{\Gw}{}u\CL^*_\gm\gz d\gg_{\gm}^\Gw=kc_\gm\gz(0)\quad\text{for all $\gz\in \BBX_\gm(\Gw)$},
  \ee
where and in the sequel the test function space 
 $$ \BBX_\gm(\Gw)=\left\{\gz\in C(\overline\Gw):\, \gg_\gm^\Gw\gz\in H_\gm(\Gw)\text{ and }\,\gr\CL^*_\gm\gz\in L^\infty(\Gw)\right\}.$$
Furthermore, if  $u$ is a nonnegative 
$\CL_\gm$-harmonic function vanishing on $\prt\Gw\setminus\{0\}$, there  exists $k\geq 0$ such that $u=k\phi^\Gw_\gm$.}
\medskip

We let $\gs_{\gm}^\Gw\in H_\mu(\Gw)$ be the unique variational solution  of  
   \bel{2-0B1}\BA {lll}
\CL_\gm u=\myfrac{\gg^\Gw_\gm}{\gr^*}\quad\text{in }\,\Gw\quad{\rm and}\quad\  u=0\quad \text{on }\, \prt\Gw,
\EA\ee
where $\rho^*(x)=\min\{\frac1{l^\Gw_\gm},\, \rho\}$. We prove that there is $c_2>1$ such that
   \bel{1-A23}\BA {lll}
    \gg^\Gw_\mu \leq  \gs_{\gm}^\Gw\leq c_2\gg^\Gw_\mu\quad{\rm in}\ \, \Omega.
\EA\ee
Note that  $\gs_{\gm}^\Gw \in C^{2}(\overline\Gw\setminus\{0\})$ is a positive classical solution of (\ref{2-0B1})
with zero Dirichlet boundary condition on $\prt \Gw\setminus \{0\}$,  i.e. $\gs_{\gm}^\Gw=0$ on $\prt\Gw\setminus\{0\}$. Moreover,  $\myfrac{\prt \gs_{\gm}^\Gw}{\prt{\bf n}}<0$ on $\prt\Gw\setminus\{0\}$. We set
$$
\eta=\myfrac{\gs_{\gm}^\Gw}{\gg^\Gw_\gm}\quad {\rm in}\ \Omega,
$$
 which satisfies 
     \bel{2-0B3}\BA {lll}
\CL^*_\gm\eta=\myfrac{1}{\gr}\quad\text{in }\,\Gw, 
\EA\ee
play a key role in the sequel. Clearly, $\eta\in C^{2}(\overline\Gw\setminus\{0\}) $ and $1\leq \eta\leq c_2$ in $\overline\Gw\setminus\{0\}$ by (\ref{1-A23}). 
We   denote by $\frak M(\Gw;\gs_{\gm}^\Gw)$ the set of Radon measures $\gn$ in $\Gw$ such that
$$
 \sup\left\{\myint{\Gw}{}\gz d|\gl|:\gz\in C_c(\Gw),\,0\leq\gz\leq\gs_{\gm}^\Gw \right\}:=\myint{\Gw}{}\gs_{\gm}^\Gw d|\gn|<+\infty.
$$
If $\gn\in \frak M_+(\Gw;\gs^\Gw_\gm)$ the measure $\gs^\Gw_\gm\gn$ is a nonnegative bounded measure in $\Gw$. Put
\bel{1-A29}
\gb^\Gw_\gm(x)=-\frac{\prt\gg_{\gm}^\Gw(x)}{\prt{\bf n_x}}=\lim_{t\to 0^+}\frac{\gg_{\gm}^\Gw(x-tn_x)}t=\lim_{t\to 0^+}\frac{\gg_{\gm}^\Gw(x-tn_x)}{\rho^*(x-tn_x))},\quad\forall\, x\in\prt\Gw\setminus\{0\}
\ee
and from (\ref{1-A23}) and (\ref{1-C1}), we have that
\bel{1-A28}\BA {lll}
c_1 |x|^{\alpha_+-1}\leq \gb^\Gw_\gm(x)\leq c_1c_2|x|^{\alpha_+-1}\quad {\rm for}\ x\in\prt\Omega\setminus\{0\}.
 \EA\ee
As a consequence, the following potential function plays an important role in defining our boundary data. Denote
\bel{1-A28-1}\BA {lll}
\gb_\gm(x)=|x|^{\ga_+-1} \quad {\rm for}\ x\in\R^N\setminus\{0\}.
 \EA\ee
 The set Radon measures $\gl$ in $\prt\Gw\setminus\{0\}$ such that 
 $$
 \sup\left\{\myint{\prt\Gw\setminus\{0\}}{}\gz d|\gl|:\gz\in C_c(\prt\Gw\setminus\{0\}),\,0\leq\gz\leq \gb_\gm \right\}:=
\myint{\prt\Gw\setminus\{0\}}{} \gb_\gm d|\gl|<+\infty
$$
is denoted by $\frak M(\prt\Gw\setminus\{0\}; \gb _\gm)$. The extension of $\gl\in \frak M_+(\prt\Gw\setminus\{0\}; \gb_\gm)$  as a measure $ \gb_\gm\gl$ in $\prt\Gw$ is given by
 $$
\myint{\prt\Gw}{}\gz d (\gb _\gm\gl)= \sup\left\{\myint{\prt\Gw}{}\gu  \gb _\gm\, d\gl:\gu\in C_c(\prt\Gw\setminus\{0\}),\,0\leq\gu\leq\gz \right\}
\quad\text{for all }\,\gz\in C_c(\prt\Gw)\,,\;\gz\geq 0
$$
and $ \gb_\gm\gl= \gb_\gm\gl_+- \gb_\gm\gl_-$ if $\gl$ is a signed measure in $\frak M(\prt\Gw\setminus\{0\}; \gb _\gm)$, and this defines the set $\frak M(\prt\Gw; \gb _\gm)$ of all such extensions. The Dirac mass at $0$ does not belong  to $\frak M(\prt\Gw; \gb _\gm)$, but it is the limit of sequences of measures in this space in the same way as it is a limit of measures in $\frak M_+(\prt\Gw\setminus\{0\}; \gb _\gm)$.  In the next result we prove the existence and uniqueness of a solution to
\bel{1-B1}\left\{\BA {lll}
\CL_\gm u=\gn\quad&\text{in }\ \Gw,\\[0.5mm]
\phantom{\CL_\gm}
u=\gl+k\gd_0\quad&\text{on }\,\prt\Gw.
\EA\right.\ee
Thanks to (\ref{1-A12}) the Green kernel $G^\Gw_\gm$ is easily constructible. If  $\gn\in \frak M_+(\Gw;\gs_{\gm}^\Gw)$ and $\gl\in \frak M(\prt\Gw; \gb _\gm)$ the following expressions are well defined
 $$
\BBK^\Gw_\gm[\gl](x)=\myint{\prt\Gw}{}K^\Gw_\gm(x,y)  d\gl(y)\quad\text{and }\;\BBG^\Gw_\gm[\gn](x)=\myint{\Gw}{}G^\Gw_\gm(x,y)  d\gn(y).
$$
Our main existence result is the following.\medskip

\nind{\bf Theorem D} {\it Let $\Gw$ be a $C^2$ bounded domain such that $0\in\prt\Gw$ satisfying $(\CC$-$1)$ and $\gm\geq\gm_1$. If $\gn\in \frak M_+(\Gw;\gs_{\gm}^\Gw)$, $\gl\in \frak M(\prt\Gw; \gb _\gm)$ and $k\in\BBR$, the function 
 \bel{1-A32}
 u=\BBG^\Gw_\gm[\gn]+\BBK^\Gw_\gm[\gl]+k\phi^\Gw_\gm:=\BBH^\Gw_\gm[(\gn,\gl,k)]
 \ee
  is the unique solution of (\ref{1-B1}) in the very weak sense that   $u\in L^1(\Gw,\gr^{-1}d\gg_\gm^\Gw)$ and
  \bel{1-Ax32}
\myint{\Gw}{}u\CL^*_\gm\gz d\gg_{\gm}^\Gw=\myint{\Gw}{}\gz d(\gg^\Gw_\gm\gn)+\myint{\prt\Gw}{}\gz d(\gb^\Gw_\gm\gl)
+kc_\gm\gz(0)
  \ee
for all $\gz\in \BBX_\gm(\Gw)$.}
\medskip

In the next result we prove that all the positive $\CL_\gm$-harmonic functions in $\Gw$ are described by formula (\ref{1-A32}) (with $\gn=0$).
\medskip

\nind{\bf Theorem E} {\it Let $\Gw$ be a $C^2$ bounded domain such that $0\in\prt\Gw$ satisfying $(\CC$-$1)$,  $\gm\geq\gm_1$ and $u$ be a nonnegative $\CL_\gm$-harmonic functions in $\Gw$. Then there exist
$\gl\in \frak M(\prt\Gw; \gb _\gm)$ and $k\geq 0$, such that 
 $$
 u=\BBK^\Gw_\gm[\gl]+k\phi^\Gw_\gm=\BBH^\Gw_\gm[(0,\gl,k)].
 $$
 The couple $(\gl,k\gd_0)$ is called the boundary trace of $u$.
 }\medskip
 
 The rest of this paper is organized as follows. In section 2, we introduce the distributional identity of $\CL_\mu$ harmonic function $\phi_\mu$ in $\R^{N}_+$.   Section 3 is devoted to build the Kato's type inequalities, to construct Poisson kernel and  related properties.   
Section 4 is addressed to classify the boundary isolated singular   $\CL_\mu$ harmonic functions in a bounded domain, i.e. {\bf Theorem A} and to show the existence and related distributional identity in a $(\CC$-$1)$ domain:  proofs of {\bf  Theorem B} and {\bf Corollary C}.   We classify the boundary trace for general  $\CL_\mu$ harmonic functions and
give the existence of $\CL_\mu$ harmonic functions with the boundary trace $(\gl,k\gd_0)$:
{\bf Theorem D} and {\bf Theorem E} respectively in Section 5. Finally, we show Estimates (\ref{1-C1}) and (\ref{1-A23}) in Appendix.

 In a forthcomming article \cite{ChVe1} we study the semilinear problem
$$
\left\{\BA {lll}
\CL_\gm u+ g(u)=0&\quad\text{in }\ \Gw,\\
\phantom{\CL_\gm + g(u)}
u=\gl&\quad\text{on }\,\prt\Gw.
\EA\right.$$




\mysection{The half-space setting}

Let
$\BBR^N_+:=\{x=(x',x_N)\in\R^{N-1}\times\R:\, x_N>0\}$, $(r,\gs)\in \BBR_+\ti \BBS^{N-1}_+$ be the spherical coordinates in $\BBR^N_+$ and $\Gd'$ is the Laplace Beltrami operator on $\BBS^{N-1}$. Then 
$$
\CL_\gm u=-\prt_{rr}u-\myfrac{N-1}{r}\prt_{r}u-\myfrac{1}{r^2}\Gd'u+\myfrac{\gm}{r^2}u.
$$
If  $u(r,\gs)=r^{\ga}\phi(\gs)$ is a (separable) solution of $\CL_\gm u=0$ vanishing on $\prt\BBR^N_+\setminus\{0\}$, then $\phi$ satisfies
$$
\left\{\BA {lll}-\Gd'\phi=\gl_k\phi&\quad\text{in }\, \BBS^{N-1}_+:=\BBS^{N-1}\cap\BBR^N_+,\\[1mm]
\phantom{-\Gd'}\phi=0&\quad\text{in }\, \prt \BBS^{N-1}_+\approx \BBS^{N-2},
\EA\right.$$
where $\gl_k$ a constant which necessarily belongs to the 
spectrum 
$$\gs_{\BBS^{N-1}_+}(-\Gd')=\{\gl_k=k(N+k-2):k\in\BBN^*\},$$ 
and $\ga=\ga_{k\,+},\ga_{k\,-}$ is a root of 
\bel{2-A3}
\ga^2+(N-2)\ga-\gl_k-\gm=0.
\ee
The fundamental state corresponds to $k=1$, in which case  since $\gl_1=N-1$, existence of real roots of (\ref{2-A3}) necessitates 
$\gm\geq \gm_1=-\frac {N^2}{4}=\gm_1$ and we denote $\ga_{1\,+}=\ga_{+}$ and $\ga_{1\,-}=\ga_{-}$.   Note that this value is connected to the  boundary Hardy
$$
\myint{\BBR^{N}_+}{}|\nabla \phi|^2+\gm_1\myint{\BBR^{N}_+}{}\myfrac{\phi^2}{|x|^2}dx\geq 0\quad\text{for all }\gf\in C^\infty_0(S^{N-1}_+). 
$$
If this condition is fulfilled, the two roots $\ga_+$ and $\ga_-$ corresponding to $k=1$ and $\gl_1$ are 
\bel{2-A5}
\ga_+=\frac{2-N}{2}+\sqrt{\gm-\gm_1}\quad\text{and }\ga_-=\frac{2-N}{2}-\sqrt{\gm-\gm_1}.
\ee
The corresponding positive separable solutions $\gg_\gm$ and $\gf_\gm$ of $\CL_\gm u=0$ vanishing on $\prt\BBR^N_+\setminus\{0\}$ are defined by (\ref{1-A13}). We set $d\gg_\gm(x)=\gg_\gm(x) dx$ and define the operator $\CL^*_\gm$ by (\ref{1-A17}).
\bprop{weak1} The function $\gf_\gm$ belongs to $L^1_{loc}(\BBR^N_+,\gr^{-1}d\gg_\gm)$. It satisfies 
\bel{2-A6}
\myint{\BBR^N_+}{}\gf_\gm\CL^*_\gm\gz d\gg_\gm(x)=c_\gm\gz(0) \quad\text{for all $\gz\in  \BBX_\gm(\BBR^N_+)$,}
\ee
where
 \bel{2-A7}\BA {lll}
c_\gm=\left\{\BA{lll}2\sqrt{\gm-\gm_1}\myint{\BBS^{N-1}_+}{}\psi_1^2dS&\quad\text{if }\;\gm>\gm_1,\\[2mm]
\myint{\BBS^{N-1}_+}{}\psi_1^2dS&\quad\text{if }\;\gm=\gm_1
\EA\right.
\EA\ee
and 
 $  \BBX_\gm(\BBR^N_+)=\left\{\gz\in C_c(\overline{\BBR^N_+}):\, \gr\CL^*_\gm\gz\in L^{\infty}(\BBR^N_+)\right\}$. 
\es
\noindent{\bf Proof.} Let $\gz\in  \BBX_\gm(\BBR^N_+)$, $\ge>0$ and set $B_\ge^+=B_\ge(0)\cap \BBR^N_+$, $\left(B_\ge^+\right)^c=B^c_\ge(0)\cap \BBR^N_+$ and $\Gg^+_\ge=\prt B_\ge(0)\cap \overline \BBR^N_+$
$$\BA {lll}
0=\myint{\left(B_\ge^+\right)^c}{}\gz\gg_\gm\CL_\gm\gf_\gm dx\\[4mm]
\phantom{0}
=\myint{\left(B_\ge^+\right)^c}{}\gf_\gm\CL_\gm^*\gz d\gg_\gm(x)+\myint{\Gg^+_\ge}{}
\left(-\myfrac{\prt\gf_\gm}{\prt{\bf n}}\gz\gg_\gm+\left(\gg_\gm\myfrac{\prt\gz}{\prt{\bf n}}+\gz\myfrac{\prt\gg_\gm}{\prt{\bf n}}\right)\gf_\gm\right) dS\\[4mm]
\phantom{0}
=\myint{\left(B_\ge^+\right)^c}{}\gf_\gm\CL_\gm^*\gz d\gg_\gm(x)+\gz(0)A(\ge)
\EA$$
and ${\bf n}=-\ge^{-1}x$ on $\Gg^+_\ge$, and
\bel{Ae}
A(\ge)=\left\{\BA {lll}-2\sqrt{\gm-\gm_1}\myint{\BBS_+^{N-1}}{}\gf^2_1dS+O(\ge)&\quad\text{if }\;\gm>\gm_1,\\[4mm]
-\myint{\BBS_+^{N-1}}{}\gf^2_1dS+O(\ge)&\quad\text{if }\;\gm=\gm_1,
\EA\right.
\ee
which implies (\ref{2-A6})-(\ref{2-A7}).\qeda

\mysection{The Poisson kernel}

In this section we assume that $\Gw$ is a bounded $C^2$ domain included in $B_{1}$ (which can always be assumed by scaling) and $0\in\prt\Gw$. 
We start with the following identity of commutation valid for all $\gl\in\mathfrak M(\prt\Gw;\gb_\gm^\Gw)$ and $\gz\in C^{1,1}(\overline\Gw)$
 \bel{2-0B0'}\BA {lll}
-\myint{\prt\Gw}{}\myfrac{\prt (\gz\gg^\Gw_\gm)}{\prt{\bf n}}d\gl=\myint{\prt\Gw}{}\gz d(\gb_\gm^\Gw \gl)\quad\text{for all }\,\gz\in C^{1,1}(\overline\Gw),
\EA\ee
where $\gb_\gm^\Gw$ is defined in (\ref{1-A29}).

The following inequality extends the classical Kato inequality to our framework.
\blemma{kato} Assume $N\geq 3$ and $\gm\geq\gm_1$ or $N=2$ and  $\gm>\gm_1$. Then for any $f\in L^1(\Gw,\gs_\gm^\Gw dx)$, $h\in L^1(\prt\Gw,\gb_\gm^\Gw dx)$ there exists a unique weak solution $u$ to
\bel{k1}\left\{\BA {lll}
\CL_\gm u=f\quad&\text{in }\ \Gw,\\[0.5mm]
\phantom{\CL_\gm }
u=h\quad&\text{on }\,\prt\Gw
\EA\right.
\ee
in the sense that  
\bel{1-Ax32-kato}
\myint{\Gw}{}u\CL^*_\gm\gz d\gg_{\gm}^\Gw=\myint{\Gw}{}f\gz d\gg^\Gw_\gm +\myint{\prt\Gw}{}h\gz d \gb^\Gw_\gm  \quad \text{for all } \gz\in \BBX_\gm(\Gw).
  \ee

Furthermore, for any $\gz\in\BBX_\gm(\Gw)$, $\gz\geq 0$, there holds
\bel{k2}
\myint{\Gw}{}|u|\CL^*_\gm\gz d\gg^\Gw_\gm\leq \myint{\Gw}{}\text{sgn}(u)f\gz d\gg^\Gw_\gm+\myint{\prt\Gw}{}|h|\gz d\gb^\Gw_\gm
\ee
and
\bel{k3}
\myint{\Gw}{}u^+\CL^*_\gm\gz d\gg^\Gw_\gm\leq \myint{\Gw}{}\text{sgn}_+(u)f\gz d\gg^\Gw_\gm+\myint{\prt\Gw}{}h^+\gz d\gb^\Gw_\gm.
\ee
\es
\noindent{\bf Proof.} {\it Uniqueness}. Assume that  $u$ is a weak solution of (\ref{k1}) with $f=h=0$. Then for any $\gz\in\BBX_\gm(\Gw)$ there holds
$$
\myint{\Gw}{}u\CL^*_\gm\gz d\gg^\Gw_\gm=0.
$$
Let $\gf\in C^\infty_c(\Gw)$, and $\gu\in H_\mu(\Gw)$ be the variational solution of
$$\CL_\gm\gu=\frac{\gg^\Gw_\gm}{\gr}\gf\quad{\rm and}\quad u\in H_\mu(\Gw).$$
  Then $\gu\in C^\infty(\Gw)$, $|\gu|\leq \norm\gf_{L^\infty}\sigma_\mu^\Gw$; the equation is satisfied everywhere and in the sense of distributions in $\Gw$. Clearly $w=(\gg^\Gw_\gm)^{-1}\gu$ belongs to $C^\infty(\Gw)$ and satisfies 
 $$\CL^*_\gm w=\frac{1}{\gr}\gf.$$ 
Thus,
$$
\myint{\Gw}{}\frac u\gr\gf d\gg^\Gw_\gm=0.
$$
Since $\gf$ is arbitrary, we have that $u=0$.\smallskip

\nind {\it Existence and estimates}. We proceed by approximation as in \cite[Prop. 2.1]{CQZ}. We assume that 
$\{(f_n,h_n)\}\subset C^1_0(\Gw)\ti C^1_0(\prt\Gw\setminus\{0\})$ is a sequence which converges to $(f,h)$ in 
$L^1(\Gw,\gg_\gm^\Gw dx)\ti L^1(\prt\Gw,\gb_\gm^\Gw dx)$. We set $V(x)=|x|^{-2}$, denote by $\BBK^\Gw$ the Poisson potential of $-\Delta$ in 
$\Gw$  and consider the approximate problem
\bel{k5}\left\{\BA {lll}
\CL_\gm w_n=f_n-\gm V\BBK[h_n]\quad&\text{in }\ \Gw,\\[0.8mm]
\phantom{\CL_\gm }
w_n=0\quad&\text{on }\,\prt\Gw.
\EA\right.
\ee
Near $0$, we have $V\BBK[h_n](x)=O(\frac{\gr(x)}{|x^2|})$, hence, if $N\geq 3$, $V\BBK[h_n]\in L^2(\Gw)$. If $N=2$, the function 
$x\mapsto \frac{\gr(x)}{|x|^2}$ belongs to the Lorentz space $L^{2,\infty}(\Gw)$ which is the dual of $L^{2,1}(\Gw)$. Since 
$H(\Gw)\subset L^{2,1}(\Gw)$ by (\ref{1-A12}), it follows that $H'(\Gw)\subset L^{2,\infty}(\Gw)$. Hence, by Lax-Milgram theorem there exists a unique $w_n\in H(\Gw)$ such that (\ref{k5}) holds in the variational sense. Then $u_n=w_n+\BBK[h_n]$, which has the same regularity as $w_n$, satisfies 
\bel{k6}\left\{\BA {lll}
\CL_\gm u_n=f_n\quad&\text{in }\ \Gw,\\[0.5mm]
\phantom{\CL_\gm }
u_n=h_n\quad&\text{on }\,\prt\Gw.
\EA\right.
\ee
For $\gs>0$, we set 
$$m_\gs(t)=\left\{\BA {lll}|t|-\frac\gs 2&\quad\text{if }\,|t|\geq\gs,\\[0.8mm]
\frac{t^2}{2\gs }&\quad\text{if }\,|t|<\gs.
\EA\right.
$$
The $m_\gs$ is convex, $|m_\gs'(t)|\leq 1$ and $m_\gs'(t)\to\,$sign$_0(t)$ as $\gs\downarrow 0$. Let $\gz\in C^{1,1}(\overline\Gw)$, $\gz\geq 0$. We have that
$$
\BA {lll}
\myint{\Gw}{}\left(\langle\nabla w_n,\nabla(\gz m_\gs'(u_n)\gg_\gm^\Gw)\rangle+ \gm Vw_n\gz m_\gs'(u_n)\gg_\gm^\Gw\right)dx\\[4mm]
\phantom{\langle\nabla w_n,\nabla(\gz m_\gs'(u_n)\gg_\gm^\Gw)\rangle+ \gm }
=\myint{\Gw}{}\gz m_\gs'(u_n)\gg_\gm^\Gw f_n dx
-\gm\myint{\Gw}{} V\BBK[h_n]\gz m_\gs'(u_n)\gg_\gm^\Gw dx=:\BBR(\gs).
\EA$$
By  the fact $u_n=w_n+\BBK[h_n]$,   we have that
$$
\BA {lll}
\BBR(\gs)&=\myint{\Gw}{}\langle\nabla u_n,\nabla(\gz m_\gs'(u_n)\gg_\gm^\Gw)\rangle dx-\myint{\Gw}{} \langle\nabla \BBK[h_n],\nabla(\gz m_\gs'(u_n)\gg_\gm^\Gw)\rangle dx\\[4mm]
&\phantom{-\,}
+\gm\myint{\Gw}{}Vu_n\gz m_\gs'(u_n)\gg_\gm^\Gw dx-\gm\myint{\Gw}{}V\BBK[h_n]\gz m_\gs'(u_n)\gg_\gm^\Gw dx\\[4mm]
&=\myint{\Gw}{}|\nabla u_n|^2 m_\gs''(u_n)\gz\gg_\gm^\Gw dx+\myint{\Gw}{}\langle\nabla m_\gs(u_n),\nabla(\gz\gg_\gm^\Gw)\rangle dx+\gm\myint{\Gw}{}Vu_n\gz m_\gs'(u_n)\gg_\gm^\Gw dx\\[4mm]
& \phantom{-\, }
-\gm\myint{\Gw}{}V\BBK[h_n]\gz m_\gs'(u_n)\gg_\gm^\Gw dx\\[4mm]
&\geq -\myint{\Gw}{} m_\gs(u_n)\Gd(\gz\gg_\gm^\Gw) dx+\myint{\prt\Gw}{} m_\gs(h_n)\gz\myfrac{\prt \gg_\gm^\Gw}{\prt{\bf n}} dS
+\gm\myint{\Gw}{}Vu_n\gz m_\gs'(u_n)\gg_\gm^\Gw dx\\[4mm]
&\phantom{ -\, }
-\gm\myint{\Gw}{}V\BBK[h_n]\gz m_\gs'(u_n)\gg_\gm^\Gw dx\\[4mm]
\phantom{}
&\geq \myint{\Gw}{} m_\gs(u_n)\CL^*_\gm\gz d\gg_\gm^\Gw -\myint{\prt\Gw}{} m_\gs(h_n)d\gb_\gm^\Gw
+\gm\myint{\Gw}{}V\left(u_n m_\gs'(u_n)-m_\gs(u_n)\right)\gg_\gm^\Gw \gz dx
\\[4mm]
&\phantom{-\, }
-\gm\myint{\Gw}{}V\BBK[h_n]\gz m_\gs'(u_n)\gg_\gm^\Gw dx,
\EA$$
thus,  we obtain that 
\bel{k8}\BA {lll}
\myint{\Gw}{} m_\gs(u_n)\CL^*_\gm\gz d\gg_\gm^\Gw+\gm\myint{\Gw}{}V\left(u_n m_\gs'(u_n)-m_\gs(u_n)\right) \gz d\gg_\gm^\Gw\\[4mm]\phantom{-------------}
\leq \myint{\Gw}{}\gz m_\gs'(u_n) f_n d\gg_\gm^\Gw+\myint{\prt\Gw}{} m_\gs(h_n)\gz d\gb_\gm^\Gw.
\EA\ee
Since $m_\gs$ is convex, $u_n m_\gs'(u_n)-m_\gs(u_n)\geq 0$. Hence for $\gm\geq 0$, we can let $\gs\to 0$ in (\ref{k8}) and obtain (\ref{k2}). 

For  $\gm\in[\gm_1,\,0)$, we note that 
$$0\leq u_n m_\gs'(u_n)-m_\gs(u_n)\leq \frac{|u_n|^2}{2\gs^2}\chi_{_{\{|u_n|\leq \gs\}}}$$ 
and
\bel{k9}0\leq \myint{\Gw}{}V\left(u_n m_\gs'(u_n)-m_\gs(u_n)\right) \gz d\gg_\gm^\Gw\leq 
\frac{\norm\gz_{L^\infty}}{2}\myint{\{|u_n|\leq \gs\}}{}|x|^{-1-\frac N2+\sqrt{\gm-\gm_1}}dx.
\ee
Hence if $N\geq 3$, or $N=2$ and $\gm>\gm_1=-1$, the right-hand side of (\ref{k9}) tends to $0$ as $\gs\to 0$ and then we obtain (\ref{k2}). The proof of (\ref{k3}) is similar.
 
Applying  estimate (\ref{k2}) to $u_n-u_m$, we obtain for all $\gz\in \BBX_\gm(\Gw)$, $\gz\geq 0$,
$$
\myint{\Gw}{}|u_n-u_m|\CL^*_\gm\gz d\gg^\Gw_\gm\leq \myint{\Gw}{}|f_n-f_m|\gz d\gg^\Gw_\gm+\myint{\prt\Gw}{}|h_n-h_m|\gz d\gb^\Gw_\gm,
$$
For  test function, we take $\eta$, the solution of (\ref{2-0B3}),
then
$$
\myint{\Gw}{}\myfrac{|u_n-u_m|}{\gr}d\gg^\Gw_\gm\leq \myint{\Gw}{}|f_n-f_m|d\gs^\Gw_\gm+\myint{\prt\Gw}{}|h_n-h_m|d(\eta\gb^\Gw_\gm),
$$
Therefore $\{u_n\}$ is a cauchy sequence in $L^1(\Gw,\gr^{-1}d\gg^\Gw_\gm)$ with limit $u$. Since $u_n$ satisfies (\ref{k6}), we let $n$ go to infty in 
 $$
\myint{\Gw}{}u_n\CL^*_\gm\gz\, d\gg_{\gm}^\Gw=\myint{\Gw}{}\gz \gs^\Gw_\gm f_ndx+\myint{\prt\Gw}{}\gz  \gb^\Gw_\gm h_ndS
  $$
 and obtain (\ref{1-Ax32-kato}).
\qeda

\blemma{id} Assume $\gl\in \mathfrak M(\prt\Gw;\gb_\gm)$ and $\gz\in\BBX_\gm(\Gw)$, then there holds
 $$
 \BA {lll}
\myint{\Gw}{}\left(\CL_\gm^*\gz-\myfrac{\gm}{|x|^2}\gz\right)\BBK^\Gw[\gl]d\gg^\Gw_\gm=\myint{\prt\Gw}{}\gz d(\gb_\gm^\Gw \gl).
\EA$$
\es
\noindent{\bf Proof.} Note that $d(\gb_\gm^\Gw \gl)$ is equivalent to $ d(\gb_\gm^\Gw \gl)$
by (\ref{1-A28}). By (\ref{2-0B0}) we have almost everywhere in $\Gw$,
$$\left(\CL_\gm^*\gz-\myfrac{\gm}{|x|^2}\gz\right)\BBK^\Gw[\gl]\gg^\Gw_\gm
=\left(\CL_\gm(\gg^\Gw_\gm\gz)-\myfrac{\gm}{|x|^2}\gg^\Gw_\gm\gz\right)\BBK^\Gw[\gl]=-\Gd(\gg^\Gw_\gm\gz)\BBK^\Gw[\gl].
$$
If we assume that $\gl$ vanishes in a neighborood of $0$ we derive from (\ref{2-0B0'})
$$-\myint{\Gw}{}\Gd(\gg^\Gw_\gm\gz)\BBK^\Gw[\gl] dx=-\myint{\prt\Gw}{}\myfrac{\prt (\gz\gg^\Gw_\gm)}{\prt{\bf n}}d\gl=\myint{\prt\Gw}{}\gz d(\gb_\gm^\Gw \gl).
$$
Since $\gg_\gm^\Gw\left(\CL_\gm^*\gz-\myfrac{\gm}{|x|^2}\gz\right)$ is bounded, we obtain the result first if $\gl$ is nonnegative
 by considering the sequence  $\{\chi_{_{B^c_\ge}}\gl\}$ and letting $\ge\to 0$, and then for any $\gl=\gl^+-\gl^-$.
\qeda\medskip

We observe also that the existence of the Green kernel follows from Lax-Milgram theorem which gives the existence of a unique variational solution in $H(\Gw)$ to 
$$\left\{\BA {lll}-\Gd u+\myfrac{\gm}{|x|^2}u=f\quad&\text{in }\ \Omega,\\[0.9mm]
\phantom{-\Gd u+\myfrac{\gm}{|x|^2}}u=0\quad&\text{on }\,\prt\Gw.
\EA\right.$$
We denote by $G^\Gw_{\gm}$ the Green kernel and by $\BBG^\Gw_{\gm}$ the corresponding Green operator. 

\subsection{Construction of the Poisson kernel when $\gm> 0$}
For the sake of completeness, we recall the construction in \cite{VeYa}. For $\ge>0$ we set $V_\ge(x)=\max\{\ge^{-2},|x|^{-2}\}$ and $V_0(x)=V(x)=|x|^{-2}$, and if 
  $\gl\in\mathfrak M(\prt\Gw)$ let $u_\ge$ be the solution  of 
 $$
 \left\{\BA {lll}
-\Gd u+\gm V_\ge u=0\quad&\text{in }\ \Gw,\\[0.5mm]
 \phantom{-\Gd u+\gm V_\ge }
u=\gl\quad&\text{on }\, \prt\Gw.
\EA\right.$$
Then
$$u_\ge(x)=\myint{\prt\Gw}{}K^\Gw_{\gm,\ge}(x,y)d\gl(y)=\BBK^\Gw_{\gm,\ge}[\gl].
$$
We obtain by the maximum principle,
$$K^\Gw_{\gm,\ge}\leq K^\Gw_{\gm',\ge'}\leq K^\Gw\quad\text{for all }\;\gm\geq \gm'\geq 0\text{ and }\;\ge'\geq\ge>0,
$$
where $K^\Gw$ is the usual Poisson kernel in $\Gw$ and there exists 
$$
K^\Gw_{\gm}(x,y)=\lim_{\ge\to 0}K^\Gw_{\gm,\ge}(x,y)\qquad\text{for all }\,(x,y)\in\Gw\ti\prt\Gw.
$$
Therefore we infer, firstly by monotone convergence if $\gl\geq 0$, and then for any  $\gl\in\mathfrak M(\prt\Gw)$, that
 \bel{2-B4}
\lim_{\ge\to 0}u_\ge(x)=u(x)=\myint{\prt\Gw}{}K^\Gw_{\gm}(x,y)d\gl(y)\qquad\text{for all }\,x\in\Gw.
\ee
Since $V$ is finite in $B^c_\ge\cap\Gw$, for any $x\in\Gw$, $K^\Gw_{\gm}(x,y)>0$ for all $y\in \prt\Gw\setminus\{0\}$.  If $G^\Gw$ is the Green kernel in $\Gw$, there holds
$$
\BA {lll}
u_\ge(x)+\gm\myint{\Gw}{}G^\Gw(x,y)V_\ge(y)u_\ge(y)dy=\myint{\prt\Gw}{}K^\Gw(x,y)d\gl(y)
\EA$$
If $\gl\geq 0$, we have by Fatou's lemma,
 \bel{2-B6}\myint{\Gw}{}G^\Gw(x,y)V(y)u(y)dy\leq\liminf_{\ge\to 0}\myint{\Gw}{}G^\Gw(x,y)V_\ge(y)u_\ge(y)dy.
\ee
Combined with (\ref{2-B4}) it yields
 $$
 \BA {lll}
u(x)+\gm\myint{\Gw}{}G^\Gw(x,y)V(y)u(y)dy\leq \myint{\prt\Gw}{}K^\Gw(x,y)d\gl(y)\qquad\text{for all }\,x\in\Gw.
\EA$$
Since the function $u+\gm\BBG[Vu]$ is nonnegative and harmonic in $\Gw$, it admits a boundary trace which is a nonnegative Radon measure $\gl^*$ and there holds
 \bel{2-B8}\BA {lll}
u(x)+\gm\myint{\Gw}{}G^\Gw(x,y)V(y)u(y)dy= \myint{\prt\Gw}{}K^\Gw(x,y)d\gl^*(y)\qquad\text{for all }\,x\in\Gw.
\EA\ee
Because of (\ref{2-B6}) $0\leq\gl^*\leq\gl$. The measure $\gl^*$ is the reduced measure associated to $\gl$. Since (\ref{2-B8}) is equivalent to 
$$u(x)=\myint{\prt\Gw}{}K^\Gw_{\gm}(x,y)d\gl^*(y),
$$
there holds
$$\myint{\prt\Gw}{}K^\Gw_{\gm}(x,y)d(\gl-\gl^*)(y)=0.
$$
This implies that $\gl=\gl^*$ in $\prt\Gw\setminus\{0\}$. With the notations of \cite{VeYa}, we recall that
$$\BA {lll}\CS ing_{_V}(\Gw):=\left\{y\in\prt\Gw:\, \exists x_0\in\Omega\ \text{s.t.}\  K^\Gw_\gm(x_0,y)=0\right\}\\[4mm]
\qquad\qquad\quad\subset Z_{_V}:=\left\{y\in \prt\Gw:\myint{\Gw}{}K_0^\Gw(x,y)V(x)\gr(x)dx=\infty\right\}.
\EA
$$
Actually, if $y\in \CS ing_{_V}(\Gw)$,  $K^\Gw_\gm(x_0,y)=0$ for any $x_0\in\Gw$ by Harnack inequality. Clearly $0\in Z_V$ and if $y\neq 0$ the integral term in the definition of $Z_V$ is finite. Hence $\CS ing_{_V}(\Gw)\subset Z_{_V}=\{0\}$. Since for any truncated cone $C_{0,\gd}\Subset\Gw$ with vertex $0$, there holds
$$\myint{C_{0,\gd}}{}V(x)\myfrac{dx}{|x-y|^{N-2}}=\infty,
$$
it follows by Ancona's result \cite[Theorem A1]{VeYa} that $0\in \CS ing_{_V}(\Gw)$. Finally
$$
K^\Gw_\gm(x,0)=0\quad\text{for all }\;x\in\Gw.
$$
\subsection{Construction of the Poisson kernel when  $\gm_1\leq\gm< 0$} For $\ge>0$ and $\gl\in C(\prt\Gw)$, $\gl\geq 0$ we denote by $w=w_{\ge,\gl}$ the variational solution in $H(\Gw)$ of
$$
\left\{\BA {lll}
-\Gd w+\gm V_\ge w=-\gm V_\ge\BBK[\gl]\quad&\text{in }\ \Gw,\\[0.5mm]
 \phantom{-\Gd w+\gm V_\ge }
w=0\quad&\text{on }\, \prt\Gw.
\EA\right.$$
Then $w_{\ge,\gl}\geq 0$ and $u=u_{\ge,\gl}:=w_\ge+\BBK[\gl]$ satisfies 
 \bel{2-C2}\left\{\BA {lll}
-\Gd u+\gm V_\ge u=0\quad&\text{in }\ \Gw,\\[0.5mm]
 \phantom{-\Gd u+\gm V_\ge }
u=\gl\quad&\text{on }\, \prt\Gw.
\EA\right.\ee
Since 
$-\Gd u_{\ge,\gl}+\gm V u_{\ge,\gl}\leq 0$, there holds from \rlemma{kato}
 $$
\myint{\Gw}{}u_{\ge,\gl}\CL^*_\gm\gz d\gg_\gm^\Gw\leq \myint{\prt\Gw}{}\gl\gz d\gb_\gm^\Gw \quad{\rm for}\ \gz \in \BBX_\mu(\Omega),\ \gz\geq0
$$
and in particular
 \bel{2-C4}\BA {lll}
\myint{\Gw}{}\myfrac{u_{\ge,\gl}}{\gr} d\gg_\gm^\Gw\leq \myint{\prt\Gw}{}\gl d(\eta \gb_\gm^\Gw).
\EA\ee
If $\ge>\ge'>0$ and $\gl'>\gl>0$ we have
$$-\myfrac{\Gd u_{\ge,\gl}}{u_{\ge,\gl}}+\myfrac{\Gd u_{\ge',\gl'}}{u_{\ge',\gl'}}=\gm \left(V_{\ge'}-V_\ge\right)\leq 0.
$$
Since
$$\BA {lll}\myint{\Gw}{}\left(-\myfrac{\Gd u_{\ge,\gl}}{u_{\ge,\gl}}+\myfrac{\Gd u_{\ge',\gl'}}{u_{\ge',\gl'}}\right)(u^2_{\ge,\gl}-u^2_{\ge',\gl'})_+ dx\\[4mm]
\qquad=\myint{\{u_{\ge,\gl}\geq u_{\ge',\gl'}\}}{}\left(\left|\nabla u_{\ge,\gl}-\myfrac{u_{\ge,\gl}}{u_{\ge',\gl'}}\nabla u_{\ge',\gl'}\right|^2+
\left|\nabla u_{\ge',\gl'}-\myfrac{u_{\ge',\gl'}}{u_{\ge,\gl}}\nabla u_{\ge,\gl}\right|^2\right)dx,
\EA$$
we deduce that the function $x\mapsto \frac{u_{\ge',\gl'}}{u_{\ge,\gl}}(x)$ is constant on the set $\{x:u_{\ge,\gl}(x)> u_{\ge',\gl'}(x)\}$. If this set is non-empty we get a contradiction since it is strictly included in $\Gw$. Therefore the mapping
$$
(\ge,\gl)\mapsto u_{\ge,\gl}
$$
is decreasing in $\ge$ and increasing in $\gl$. 
\smallskip

Next we can assume that $\gl\in \mathfrak M_+(\prt\Gw)$ vanishes in $B_\gd\cap\prt\Gw$ and that $\{\gl_n\}\subset C(\prt\Gw)$ is a sequence of functions which converge to $\gl$ in the weak sense of measures. We denote by $u_{\ge,\gl_n}$ the solution of (\ref{2-C2}) with $\gl$ replaced by $\gl_n$. Since $\gm<0$, $\ga_+<1$,
 $\gr^{-1}\gg_\gm^\Gw\sim|x|^{\ga_+-1}\geq R^{\ga_+-1}_\gw$, where $R_\gw=\max\{|z|:z\in\Gw\}$. Hence
$$
\myint{\Gw}{}u_{\ge,\gl_n}dx\leq c_8\myint{\prt\Gw}{}\gl_n d(\eta \gb_\gm^\Gw)\leq c_{9}\norm\gl_{\mathfrak M(\prt\Gw)}.
$$
Hence $u_{\ge,\gl_n}$ and $V_\ge u_{\ge,\gl_n}$ are uniformly bounded in $L^1(\Gw)$. From standard regularity estimates the sequence $\{u_{\ge,\gl_n}\}_{n\in\BBN}$ is bounded in the Lorentz spaces $L^{\frac{N}{N-1},\infty}(\Gw)$ and weakly relatively compact in $L^{1}(\Gw)$ (see e.g. \cite{GmVe}). This implies that, up to a subsequence, $u_{\ge,\gl_n}$ converges in $L^1(\Gw)$
and a.e. in $\Gw$ to a weak solution $u_{\ge,\gl}$  of 
$$
 \left\{\BA {lll}
-\Gd u+\gm V_\ge u=0\quad&\text{in }\Gw,\\[0.5mm]
 \phantom{-\Gd u+\gm V_\ge }
u=\gl\quad&\text{on }\prt\Gw,
\EA\right.$$
that is a function which satisfies 
 \bel{2-C7}\BA {lll}
\myint{\Gw}{}\left(-u_{\ge,\gl}\Gd\gz+\gm V_\ge u_{\ge,\gl}\gz\right) dx=-\myint{\prt\Gw}{}\myfrac{\prt\gz}{\prt{\bf n}}d\gl\quad\text{for all }\;\gz\in C^{1,1}_0(\overline\Gw).
\EA\ee
Furthermore (\ref{2-C4}) holds (with the same notation). For test function $\gz$ in (\ref{2-C7}), we take  $\gz=\gth_1$ be the solution of 
$$
 \left\{\BA {lll}
-\Gd\gth_1=1\quad&\text{in }\ \Gw,\\[0.5mm]
\phantom{-\Gd}
\gth_1=0\quad&\text{on }\;\prt\Gw
\EA\right.$$
Then 
 $$
\myint{\Gw}{}u_{\ge,\gl}dx =-\gm \myint{\Gw}{}V_\ge u_{\ge,\gl}\gth_1dx-\myint{\prt\Gw}{}\myfrac{\prt\gth_1}{\prt{\bf n}}d\gl.
$$
By the monotone convergence theorem we obtain that $V_\ge u_{\ge,\gl}\to V u_{\gl}$ in $L^1(\Gw,\gth_1 dx)$ by letting $\ge\to 0$ and 
$$
\myint{\Gw}{}u_{\gl}dx =-\gm \myint{\Gw}{}V u_{\gl}\gth_1dx-\myint{\prt\Gw}{}\myfrac{\prt\gth_1}{\prt{\bf n}}d\gl.
$$
Hence 
 $$
\myint{\Gw}{}\left(-\Gd\gz+\gm V\gz\right) u_{\gl}dx=-\myint{\prt\Gw}{}\myfrac{\prt\gz}{\prt{\bf n}}d\gl\quad\text{for all }\;\gz\in C^{1,1}_0(\overline\Gw).
$$
We also have 
 $$
\myint{\Gw}{}u_{\ge,\gl_n}\CL^*_\gm\gz d\gg_\gm^\Gw=\gm\myint{\Gw}{}(V-V_\ge)u_{\ge,\gl_n}\gz d\gg_\gm^\Gw+\myint{\prt\Gw}{}\gz\gl_nd\gb_\gm^\Gw
$$
for all $\gz\in \BBX_\gm(\Gw)$. Since $u_{\ge,\gl_n}$ converges in $L^1(\Gw)$ we obtain if $\gz\geq 0$,
 $$
\myint{\Gw}{}u_{\ge,\gl}\CL^*_\gm\gz d\gg_\gm^\Gw=\gm\myint{\Gw}{}(V-V_\ge)u_{\ge,\gl}\gz d\gg_\gm^\Gw+\myint{\prt\Gw}{}\gz d(\gl \gb_\gm^\Gw)\leq \myint{\prt\Gw}{}\gz d(\gl \gb_\gm^\Gw),
$$
since $\gm(V-V_\ge)\leq 0$. When $\ge\to 0$, $u_{\ge,\gl}$ increases and converges to some $u_{\gl}$ in $L^1(\Gw,\gr^{-1}d\gg_\gm^\Gw)$ which satisfies
 $$
\myint{\Gw}{}u_{\gl}\CL^*_\gm\gz d\gg_\gm^\Gw\leq \myint{\prt\Gw}{}\gz d(\gl \gb_\gm^\Gw)\quad \text{for all } \gz\in \BBX_\gm(\Gw), \ \gz\geq 0.
$$
 For $\gd>0$ denote by $\gz_\gd$ the solution of 
 $$
 \CL^*_\gm\gz_\gd=\chi_{_{\Gw_\gd}}\CL^*_\gm\gz\quad\text{where }\;\Gw_\gd=\{x\in\Gw:\gr(x)>\gd\}.
$$
As  $\gz\in\BBX_\gm(\Gw)$, $|\CL^*_\gm\gz|\leq c_{10}\gr $, hence  $\gz_\gd\in\BBX_\gm(\Gw)$, $|\gz_\gd|\leq c_{10}\eta$ and $\gz_\gd\to\gz$ when $\gd\to 0$. Furthermore, since $c_{11}|x|$ is a supersolution for $c_{11}>0$ large enough, $\eta_\gd\leq c_{11}|x|$. Hence
$$
\myint{\Gw_\gd}{}u_{\ge,\gl}\CL^*_\gm\gz d\gg_\gm^\Gw= \gm\myint{\{|x|<\ge\}}{}\myfrac{1}{|x|^2}u_{\ge,\gl}\gz_\gd d\gg_\gm^\Gw+\myint{\prt\Gw}{}\gz_\gd d(\gl \gb_\gm^\Gw).
$$
Because $|x|^{-2}u_{\ge,\gl}|\gz_\gd|\leq c_{11}\gr^{-1}u_{\ge,\gl}$ and $u_{\ge,\gl}\to u_{\gl}$ in $L^1(\Gw,\gr^{-1}d\gg_\gm^\Gw)$, we derive that 
$$\lim_{\ge\to 0}\myint{\{|x|<\ge\}}{}\myfrac{1}{|x|^2}u_{\ge,\gl}\gz_\gd d\gg_\gm^\Gw= 0
$$
which implies
 $$
\myint{\Gw_\gd}{}u_{\gl}\CL^*_\gm\gz d\gg_\gm^\Gw= \myint{\prt\Gw}{}\gz_\gd d(\gl \gb_\gm^\Gw).
$$
Letting $\gd\to 0$ we obtain by monotonicity
 \bel{2-C14}\BA {lll}
\myint{\Gw}{}u_{\gl}\CL^*_\gm\gz d\gg_\gm^\Gw= \myint{\prt\Gw}{}\gz d(\gl \gb_\gm^\Gw).
\EA\ee

Finally, if $\gl\in\mathfrak M_+(\prt\Gw,\, \beta_\mu)$ we replace it by $\gl_\gd=\chi_{_{B^c_\gd}}\gl$ and denote by $u_{\gl_\gd}$ the weak solution of 
$$
\left\{\BA {lll}
-\Gd u+\gm V u=0\quad&\text{in }\ \Gw,\\[0.5mm]
 \phantom{-\Gd u+\gm V }
u=\gl_\gd\quad&\text{on }\, \prt\Gw.
\EA\right.$$
The mapping $\gd\mapsto u_{\gl_\gd}$ is monotone. Hence, by the monotone convergence theorem $u_{\gl_\gd}$ increases and converges to some $u_{\gl}$ in $L^1(\Gw,\gr^{-1}d\gg_\gm^\Gw)$ and clearly $u_{\gl}$ satisfies  (\ref{2-C14}) for all $\gz\in\BBX_\gm(\Gw)$.
\mysection{The singular kernel}

In this section we construct the singular kernel $\phi^\Gw_{\gm}$ and prove that it satisfies estimates (\ref{1-A25})-(\ref{1-A26}) and 
it is associated to Dirac mass at $0$. 
Up to a rotation we can assume that the inward normal direction to $\prt\Gw$ at $0$ is ${\bf e}_N=(0',1)\in\R^{N-1}\times\R$. Hence the tangent hyperplane to $\prt\Gw$ at $0$ is $\prt\BBR^N_+=\BBR^{N-1}$. For $R>0$ set $B'_R=\{x'\in\BBR^{N-1}:|x'|<R\}$ and $D_R=B'_R\ti (-R,R)$. Then there exist $R>0$ and a $C^2$ function $\gth: B'_R\mapsto\BBR$ such that $\prt\Gw\cap D_R=\{x=(x',x_N): \, x_N=\gth(x')\ \text{for }x'\in B'_R\}$ and $\Gw\cap D_R=\{x=(x',x_N): \gth(x')<x_N<R\}$. Furthermore $\nabla \gth(0)=0$.  

\subsection{Classification of Boundary isolated singularities}
We characterize the positive solutions of $\CL_\gm u=0$ which vanish on $\prt\Gw\setminus\{0\}$. 

\blemma{trace} Let $\gm\geq\gm_1$ and $u\in C^2(\overline\Gw\setminus\{0\})$ be a positive solution of $\CL_\gm u=0$ in 
$\Gw$ vanishing on $\prt\Gw\setminus\{0\}$. Then there exist $a>0$ and $c_{12}>0$ such that 
\bel{2-C16}\BA {lll}
u(x)\leq c_{12}|x|^{-a-1}\gr (x)\qquad\text {for all }\; x\in \overline\Gw\setminus\{0\}.
\EA\ee
\es
\noindent{\bf Proof.}  This is a direct consequence of Boundary Harnack inequality \cite[Th. 2.7]{B-VBV}. 
\bprop{trace2} Assume that $\gm\geq\gm_1$ and $u\in C^2(\overline\Gw\setminus\{0\})$ is a positive solution of $\CL_\gm u=0$, vanishing on $\prt\Gw\setminus\{0\}$ satisfying (\ref{2-C16}) with $a\geq-\ga_-$. Then the following convergences hold in $C^1(\BBS^{N-1}_+)$:\smallskip

\nind (i) If $\gm>\gm_1$ and $a=-\ga_-$, there exists $c_{13}\geq 0$ such that 
\bel{2-C17}\BA {lll}\displaystyle
\lim_{r\to 0}\myfrac{u(r,\cdot)}{r^{\ga_-}}=c_{13}\phi\qquad\text {as }\; r\to 0.
\EA\ee
 \smallskip

\nind (ii) $\gm\geq \gm_1$ and $a>-\ga_-$ there exist $\gt>a+\ga_-$ depending on $a$ and $\gm$, and $c_{14}\geq 0$ such that 
$$
u(x)\leq c_{14}|x|^{-a-1+\gt}\gr (x)\quad\text{for all } x\in\overline\Omega\setminus \{0\}.
$$
\es 
\noindent{\bf Proof.}  {\it Step 1. Straightening the boundary}. We define the function $\Gth=(\Gth_1,...,\Gth_N)$ on $D_R$ by $y_j=\Gth_j(x)=x_j$ if $1\leq j\leq N-1$ 
and $y_N=\Gth_N(x)=x_N-\gth(x')$. Since $D\Gth(0)=Id$ we can assume that $\Gth$ is a diffeomorphism from $D_R$ onto $\Gth(D_R)$. We set
\bel{2-C19}\BA {lll}\displaystyle
u(x)=\tilde u(y)\qquad\text{for all }x\in D^+_R=B'_R\ti[0,R).
\EA\ee
Then
 \bel{2-C20}
\BA {lll}\displaystyle
u_{x_jx_j}&=\tilde u_{y_jy_j}-2\gth_{x_j}\tilde u_{y_jy_N}-\gth_{x_j,x_j}\tilde u_{y_N}+\gth^2_{x_j}\tilde u_{y_Ny_N}\quad\text{for }1\leq j\leq N-1,\\
\!\!u_{x_Nx_N}&=\tilde u_{y_Ny_N}
\EA \ee
and  
 \bel{2-C21}\BA {lll}\displaystyle
\Gd\tilde u+|\nabla \gth|^2\tilde u_{y_Ny_N}-2\langle\nabla \gth,\nabla\tilde u_{y_N}\rangle-\tilde u_{y_N}\Gd \gth-\myfrac{\gm}{|\Gth^{-1}(y)|^2}\tilde u=0.
 \EA\ee
We use here the spherical coordinates $(r,\gs)$ in the variable $y$ and we recall that $\Gd'$ is the Laplace-Beltrami operator on $\BBS^{N-1}$ and $\nabla'$ is the tangential gradient on $\BBS^{N-1}$ identified with the covariant derivative via the isometric imbedding $\BBS^{N-1}\subset_{\!\!>}\BBR^N$ which enables the formula 
$$\nabla \tilde u (y)=\left(\tilde u_r{\bf n}+\myfrac{1}{r}\nabla'\tilde u\right)(r,\gs)\quad\text{with }\;{\bf n}=|y|^{-1}y.
$$ 
After a lengthy computation the details of which can be found in 
\cite [P\,298-300]{GmVe} we obtain
$$
\BA {lll}\displaystyle
r^2\tilde u_{rr} \left[1-2\gth_r\langle {\bf n},{\bf e}_N\rangle+|\nabla\gth|^2(\langle {\bf n},{\bf e}_N\rangle)^2\right]\\[2mm]
\;+r\tilde u_{r}
\left[N-1-r\langle {\bf n},{\bf e}_N\rangle\Gd\gth 
 +r|\nabla\gth|^2\left(\langle\nabla'(\langle{\bf n},{\bf e}_N\rangle),{\bf e}_N\rangle-2\langle\nabla'\gth,\nabla'(\langle {\bf n},{\bf e}_N\rangle)\rangle\right)\right]\\[2mm]
\;+\langle\nabla'\tilde u,{\bf e}_N\rangle\left[-r\Gd\gth+2\gth_r-|\nabla\gth|^2\langle{\bf n},{\bf e}_N\rangle\right]
+r\langle\nabla'\tilde u_r,{\bf e}_N\rangle\left[2\gth_r+2|\nabla\gth|^2\langle{\bf n},{\bf e}_N\rangle\right]
\\[2mm]
\;-2\langle\nabla'\tilde u,\nabla'\gth\rangle\langle{\bf n},{\bf e}_N\rangle
+\left\langle\nabla'(\langle\nabla'\tilde u,{\bf e}_N\rangle),|\nabla\gth|^2{\bf e}_N-2r^{-1}\nabla'\gth\right\rangle
+\Gd'\tilde u-\myfrac{\gm}{|\Gth^{-1}(y)|^2}=0.\EA$$
Next we set
$$
\tilde u(r,\gs)=r^{-a}v(t,\gs)\quad\text{with }\; t=\ln r,
$$
and we assume that 
\bel{2-C23x}\BA {lll}\displaystyle
a\neq \frac{N-2}{2}.
\EA\ee
We notice that 
$$\BA {lll}\displaystyle r^2=\sum_{j=1}^Ny_j^2=\sum_{j=1}^{N-1}x_j^2+(x_N-\gth (x'))^2=|x^2|-2x_N\gth(x')=|x|^2(1+O(r))\;\text{ as }\;r\to 0\\\phantom{\displaystyle r^2=\sum_{j=1}^Ny_j^2=\sum_{j=1}^{N-1}x_j^2+(x_N-\gth (x'))^2}
=|x|^2(1+O(e^t))\;\text{ as }\;t\to-\infty.
\EA$$
By a straightforward computation we find that $v$ satisfies the following asymptotically autonomous equation in $(-\infty,r_0]\ti \BBS^{N-1}_+$
\bel{2-C24}\BA {lll}\displaystyle
(1+\ge_1(t,\cdot))v_{tt}+\left(N-2-2a+\ge_2(t,\cdot)\right) v_t+\left(a(a+2-N)-\gm+\ge_3(t,\cdot)\right)v\\[2mm]
\phantom{-----}
+\Gd' v+\langle \nabla' v,\ge_4(t,\cdot)\rangle+\langle \nabla' v_t,\ge_5(t,\cdot)\rangle+\langle \nabla'(\langle\nabla' v,{\bf e}_N\rangle),\ge_6(t,\cdot)\rangle
=0,
\EA\ee
where the $\ge_j$ satisfies
\bel{2-C25}\BA {lll}\displaystyle
|\ge_j(t,\cdot)| +|\prt_t\ge_j(t,\cdot)|+|\nabla'\ge_{j}(t,\cdot)|\leq c_{15}e^t.
\EA\ee
This is due to the fact that $|\gth(x')|=O(|x'|^2)$ near $0$.\smallskip

\nind{\it Step 2. The convergence process}. Since $v$ is bounded in $(-\infty,r_0]\ti \BBS^{N-1}_+$ and vanishes on $(-\infty,r_0]\ti \prt \BBS^{N-1}_+$ and all the coefficients are continuous functions, we obtain that
$v$ is bounded in $W^{2,q}([T-1,T+1]\ti \BBS^{N-1}_+)$  independently of $T\leq r_0-2$, for any $q<\infty$. Hence $v$ is bounded in any 
$C^{1,\gt}([T-1,T+1]\ti \overline{\BBS^{N-1}_+})$ for any $\gt\in [0,1)$. Differentiating the equation and using the standard elliptic equations regularity, we obtain that $v$ is bounded in $W^{3,q}([T-1,T+1]\ti \BBS^{N-1}_+)$ and in $C^{2,\gt}([T-1,T+1]\ti \overline{\BBS^{N-1}_+})$. 
We consider the negative trajectory of $v$ in $C_0^1(\overline{\BBS^{N-1}_+})$ defined by
$$
 \CT_-(v)=\bigcup_{t\leq r_0-1}\{v(t,.)\}.
 $$
 By the previous estimates and the Arzela-Ascoli theorem, it is a relatively compact subset of  $C_0^1(\overline{\BBS^{N-1}_+})$, hence its 
limit set at $-\infty$ (or alpha-limit set), denoting  $A(\CT_-(v))$, is a non-empty connected compact subset of $C_0^1(\overline{\BBS^{N-1}_+})$. Multiplying (\ref{2-C24}) by $v_t$ and integrating on $\BBS^{N-1}_+$ yields
 \bel{2-C26}\BA {lll}\displaystyle
\myint{\BBS^{N-1}_+}{} \left(N-2-2a+  \ge_2-\myfrac{1}{2}\prt_t\ge_{1}\right)v_t^2 dS-\myfrac{1}{2}\myint{\BBS^{N-1}_+}{}\prt_t\ge_{3} v^2 dS\\[4mm]
\quad =\myfrac{d}{dt}\left[\myint{\BBS^{N-1}_+}{}\left(\myfrac{1}{2}|\nabla v|^2-\myfrac{1}{2}\left[a(a+2-N)-\mu+\ge_3\right]v^2-\myfrac{1}{2}(1+\ge_1)v_t^2 \right) dS\right]\\[4mm]
\qquad\  -\myint{\BBS^{N-1}_+}{}\left(\langle\nabla' v,\ge_4\rangle+\langle\nabla' v_t,\ge_5\rangle+
\langle\nabla'(\langle\nabla' v,{\bf e}_N\rangle),\ge_6\rangle\right)v_t^2 dS.
 \EA\ee
  
 Next we integrate over $(-\infty,r_2)$ for some $r_2$ large enough so that 
 $$\left|N-2-2a+\ge_2-\myfrac{1}{2}\prt_t\ge_{1}\right|\geq \myfrac{1}{2}\left|N-2-2a\right|>0,
 $$
 here we use the crucial assumption (\ref{2-C23x}). Since all the terms on the right-hand side of (\ref{2-C26}) are integrable on 
 $(-\infty,r_2)$ because of (\ref{2-C25}) and the bounds on $v$, we obtain that 
   \bel{2-C28}\BA {lll}\displaystyle
\myint{-\infty}{r_2}\myint{\BBS^{N-1}_+}{}v_t^2 dS<\infty.
 \EA\ee
 Differentiating (\ref{2-C24}) with respect to t and using the estimates on $v$ and the $\ge_j$ we obtain (see \cite[p. 302]{GmVe} for a similar calculation)
    \bel{2-C29}\BA {lll}\displaystyle
\myint{-\infty}{r_2}\myint{\BBS^{N-1}_+}{}v_{tt}^2 dS<\infty.
 \EA\ee
Because $v_t$ and $v_{tt}$ are uniformly continuous on $(-\infty,r_1]$, we infer from (\ref{2-C28}) and (\ref{2-C29}) 
   $$
\lim_{t\to-\infty}\left(\norm{v_{t}(t,.)}_{L^2(\BBS^{N-1}_+)}+\norm{v_{tt}(t,.)}_{L^2(\BBS^{N-1}_+)}\right)=0.
$$
 Therefore the set $A(\CT_-(v))$ is a compact connected subset of the set of nonnegative solutions of 
   $$
   \left\{\BA {lll}\displaystyle
\Gd'\gw+\left(a(a+2-N)-\gm\right)\gw&=0\qquad\text{in }\ \BBS^{N-1}_+,\\[1mm]
\phantom{\left(N-2-2a\right)\gw+\Gd'---}
\gw&=0\qquad\text{on }\,\prt \BBS^{N-1}_+.
 \EA\right.$$
 
 \nind{\it Step 3. The case $a(a+2-N)-\gm= N-1$}. The set $A(\CT_-(v))$ is a subset of $\ker (-\Gd'-(N-1))I_d$ in $H^{1}_0(\BBS^{N-1}_+)$ and more precisely $A(\CT_-(v))=\{m\psi_1:m\in I^*\}$ where $I^*$ is a compact interval of $[0,\infty)$. 
 We set
 $$X(t)=\myint{\BBS^{N-1}_+}{}v(t,.)\psi_1dS.
 $$
 Then $X$ satisfies
       \bel{2-C32}\BA {lll}\displaystyle
X''(t)+(N-2-2a)X'(t)+F(t)=0,
 \EA\ee
 where
       $$
       \BA {lll}\displaystyle
F(t)=\myint{\BBS^{N-1}_+}{}\left[\ge_1(t,.)v_{tt}+\ge_2(t,.)v_t+\ge_3(t,.)v+\langle \nabla' v,\ge_4(t,.)\rangle\right.\\[2mm]
\phantom{------------}
\left.+\langle \nabla' v_t,\ge_5(t,.)\rangle+\langle \nabla'(\langle\nabla' v,{\bf e}_N\rangle),\ge_6(t,.)\rangle
\right]\psi_1dS.
 \EA$$
 Then $|F(t)|\leq c_{16}e^t$. We consider a sequence $\{t_n\}$ converging to $-\infty$ and $c^*\in I^*$ such that $X(t_n)\to c^*$. Since $X'(t)$ and $X''(t)$ converges to $0$ as $t\to-\infty$, we integrate (\ref{2-C32}) on $(t_n,t)$ and let $n\to\infty$. Then we get
 $$X'(t)+(N-2-2a)(X(t)-c^*)+O(e^t)=0.
 $$
 Letting $t\to-\infty$ yields $X(t)\to c^*$. Hence we have proved that 
       $$
\lim_{t\to-\infty}v(t,.)=c^*\psi_1\quad\text{in }\; C^1(\overline{\BBS^{N-1}_+}).
 $$
  
 \nind{\it Step 4. The case $a(a+2-N)-\gm\neq N-1$}.
Clearly $A(\CT_-(v))=\{0\}$ and 
       \bel{2-C34}\BA {lll}\displaystyle
\lim_{t\to -\infty}v(t,.)=0\qquad\text{in } C^{1} (\overline{\BBS^{N-1}_+}).
 \EA\ee
Furthermore, since we have assumed $a\geq -\ga_-$, there holds actually $a> -\ga_-$. We recall that $\gl_k$ is the k-th eigenvalue of 
$-\Gd'$ in $H^1_0(\BBS^{N-1}_+)$ and put 
$$\displaystyle H_k=\ker(-\Gd'-\gl_kId)={\rm span} \langle \phi_{k,1},\phi_{k,2},...,\phi_{k,j_k}\rangle\; \text{ and }\; H^1_0(\BBS^{N-1}_+)=\overset{\infty}{\underset{k=1}{\oplus}} H_k.$$
We denote 
  $$
P_k(x)=x^2+(N-2)x-\gm-\gl_k.
$$
 Then $P_1(\ga_-)=0$ and $P_k(\ga_-)=\gl_1-\gl_k<0$ for $k\geq 2$. Since $a(a+2-N)-\gm\neq N-1$ by assumption,  we define a partition of 
 $\BBN^*$ by setting 
 $$N_1:=\{k\in \BBN^*:a(a+2-N)-\gm-\gl_k\geq 0\},\;\text{ }\;N_2:=\{k\in \BBN^*:a(a+2-N)-\gm-\gl_k<0\}
 $$
 and
 $$W_1={\underset{k\in N_1}{\oplus}} H_k\;\text{ and }W_2={\underset{k\in N_2}{\oplus}} H_k.
 $$
 Then
        \bel{2-C36}\BA {lll}\displaystyle
-\myint{\BBS^{N-1}_+}{}\phi\Gd' \phi dS\geq \gg\myint{\BBS^{N-1}_+}{}\phi^2dS\quad\text{for all }\;\phi\in W_2,
 \EA\ee
 where 
 $$\gg=\gm+\gl_{k_2}-a(a+2-N)>0\;\text{ with }\;k_2=\inf N_2.
 $$
We denote by $P_j$ the orthognal projector onto $W_j$ in $H^1_0(\BBS^{N-1}_+)$ and set $v=P_1v+P_2v=v_1+v_2$. Then the projection of 
(\ref{2-C24}) on to $W_2$ is
$$
(v_2)_{tt}+\left(N-2-2a\right) (v_2)_{t}+\left(a(a+2-N)-\gm\right)v_2+\Gd' v_2=F_2(t,.),
$$
where $F_2$ satisfies the same estimates (\ref{2-C25}) as $\ge_j$. 
Then, using  (\ref{2-C25}) and (\ref{2-C36})
 $$\myint{\BBS^{N-1}_+}{}(v_{2})_{tt}v_2dS+(N-2-2a)\myint{\BBS^{N-1}_+}{}(v_{2})_{t}v_2dS-\gg\myint{\BBS^{N-1}_+}{}v^2_2dS\geq -c_{17}e^t\left(\myint{\BBS^{N-1}_+}{}v^2_2dS\right)^{\frac12}.
 $$
 Put $Y(t)=\norm{v_2(t,.)}_{L^2(\BBS^{N-1}_+)}$, because
 $$\myint{\BBS^{N-1}_+}{}(v_2)_{t}v_2dS=Y'(t)Y(t) \;\text{ and }\; \myint{\BBS^{N-1}_+}{}(v_2)_{tt}v_2dS\geq Y''(t)Y(t),
 $$
 we obtain the following differential inequality
      $$
Y''+(N-2-2a)Y'-\gg Y\geq -c_{17}e^t\;\text{ in }\;\CD'(-\infty,r_2) .
$$
The characteristic roots of the equation $y''+(N-2-2a)y'-\gg y=0$ are 
         \bel{2-C37x}\BA {lll}a_{k_2,-}=a+\myfrac{1}{2}\left(2-N-\sqrt{4\gm+4\gl_{k_2}+(N-2)^2}\right)=\ga_{k_2,-}+a<0\\[4mm]
a_{k_2,+}=a+\myfrac{1}{2}\left(2-N+\sqrt{4\gm+4\gl_{k_2}+(N-2)^2}\right)=\ga_{k_2,+}+a>0.
\EA \ee
where the $\ga_{k_2,\pm}$ are the roots of equations (\ref{2-A3}) with $k=k_2$. The solutions of 
$$
z''+(N-2-2a)z'-\gg z= -c_{17}e^t\;\text{ in }\;\CD'(-\infty,r_2) .
$$
 endow the form $z(t)=Ae^{ta_{k_2,-}}+Be^{ta_{k_2,+}}+c_{18}e^{t}$ if $a_{k_2,+}\neq 1$ or 
 $z(t)=Ae^{ta_{k_2,-}}+Be^{t}+Cte^{t}$ if $a_{k_2,+}= 1$, for some explicit constant $c_{18}$ depending on $c_{17}$ and the coefficients in the equation.  Since $Y(t)\to 0$ when $t\to-\infty$ by (\ref{2-C34}), it follows from the maximum principle that 
$$
Y(t)\leq c_{19}e^{ta_{k_2\,+}}+c_{18}e^{t}\;\text{ if }\; a_{k_2\,+}\neq 1,\;\text{ or }\; Y(t)\leq c_{20}|t|e^{t}\;\text{ if }\; a_{k_2\,+}= 1\quad\text{for }\,t\leq r_2.
$$
Then using standard elliptic equations a priori estimates, initialy in $L^2(\BBS^{N-1}_+)$, then in $L^p(\BBS^{N-1}_+)$ and finally in $C^{\gt}(\overline{\BBS^{N-1}_+})$, we obtain that for $t\leq r_3$,  
 \bel{2-C39}\displaystyle
\norm{v_2(t,.)}_{C^1(\BBS^{N-1}_+)}\leq \left\{\BA {lll}c_{21}e^{ta_{k_2,+}}+c_{22}e^{t}\qquad&\text{ if }\; a_{k_2,+}\neq 1,\\[2mm]
c_{23}|t|e^{t}\qquad&\text{ if }\;a_{k_2,+}= 1,\EA\right.
\ee 
 where $r_3\leq r_2-1$. \\
 For the components in $W_1$ we have
  \bel{2-C39x}\displaystyle v_1(t,\cdot)=\sum_{k\in N_1}\sum_{1\leq j\leq j_k}w_{k,j}(t)\phi_{k,j}(\cdot),
\ee
 where the $\phi_{k,j}$ form an orthoromal basis of $H_k$. Then
  \bel{2-C40}\BA {lll}\displaystyle
w''_{k,j}+(N-2-2a)w'_{k,j} +\left(a(a+2-N)-\gm-\gl_k\right)w_{k,j}=F_{k,j}(t)
 \EA\ee 
 The characteristic roots of equation $z''+(N-2-2a)z' +\left(a(a+2-N)-\gm-\gl_k\right)z=0$ are given in (\ref{2-C37x}) with a general $k$,
$a_{k\,-}=a+\ga_{k\,-}$ and $a_{k\,+}=a+\ga_{k\,+}$ where $\ga_{k\,\pm}$ are the roots of (\ref{2-A3}).
They have same sign (including $0$) since $a(a+2-N)-\gm-\gl_k\geq 0$, furthermore, their sum is positive since 
$N - 2-2a<0$, as a consequence of  $a>-\ga_-$.
By standard calculation the solution of (\ref{2-C40}) has the form
  \bel{2-C41}\BA {lll}\displaystyle
w_{k,j}(t)=m_1e^{ ta_{k\,+}}+m_2e^{ ta_{k\,-}}-\myint{t}{0}\myfrac{e^{ (t-s)a_{k\,+}}-e^{ (t-s)a_{k\,-}}}{a_{k\,+}-a_{k\,-}}F_{k,j}(s)ds.
 \EA\ee
 Since $|F_{k,j}(s)|\leq c_{24}e^s$ there holds
    \bel{2-C42}\displaystyle
\left|\myint{t}{0}\myfrac{e^{ (t-s)a_{k\,+}}-e^{ (t-s)a_{k\,-}}}{a_{k\,+}-a_{k\,-}}F_{k,j}(s)ds\right|
\leq c_{25}\left\{\BA {lll}|t|e^t&\quad\text{ if }a_{k\,-}=1\\
\max\{e^t,e^{ta_{k\,-}}\}&\quad\text{ if }a_{k\,-}\neq 1
 \EA\right.\ee
 In particular, if $k_1=\max N_1$, then $a_{k_1\,\pm}=\min\{a_{k\,\pm}:k\in N_1\}$. \\
 We assume first that $a_{k_1\,-}>0$. Combining this fact with (\ref{2-C39x}) and 
 (\ref{2-C42}) we obtain
     \bel{2-C43}\displaystyle
\norm{v_1(t,.)}_{L^\infty(\BBS^{N-1}_+)}
\leq c_{26}\left\{\BA {lll}|t|e^t&\quad\text{ if }a_{k_1\,-}=1\\
\max\{e^t,e^{ ta_{k_1\,-}}\}&\quad\text{ if }a_{k_1\,-}\neq 1
 \EA\right.\ee
 Furthermore, because of the explicit formulation and (\ref{2-C25}), the left-hand side of (\ref{2-C42}) can be replaced by 
$ \norm{v_1(t,.)}_{C^1(\overline{\BBS^{N-1}_+})}$. Combining (\ref{2-C39}) and (\ref{2-C42}) we obtain the result since $v(t,.)=0$ on 
$(-\infty,r_1)\ti\prt \BBS^{N-1}_+$. \\
Next we suppose that $a_{k_1\,-}=0$. Then for $k=k_1$, (\ref{2-C41}) endows the form 
  \bel{2-C44}\BA {lll}\displaystyle
w_{k_1,j}(t)=m_1e^{ ta_{k_1\,+}}+m_2-\myfrac{1}{a_{k_1\,+}}\myint{t}{0}(e^{ (t-s)a_{k_1\,+}}-1)F_{k_1,j}(s)ds.
 \EA\ee
This implies that 
$$w_{k_1,j}(t)\to m_2+\myfrac{1}{a_{k_1\,+}}\myint{-\infty}{0}F_{k_1,j}(s)ds:=A_{k_1,j}\quad\text{as }\;t\to\infty.
$$
If $A_{k_1,j}\neq 0$ it would imply that $\displaystyle \sum_{j=1}^{j_k}A_{k_1,j}\phi_{k_1,j}$ is a nonzero eigenfunction of order $k_1>1$, hence it changes sign and it would imply that $v$ changes sign at $-\infty$ (notice that all the other terms $w_{k,j}(t)$ tends to $0$ exponentially because of (\ref{2-C41})-(\ref{2-C42})). Hence $A_{k_1,j}= 0$ and (\ref{2-C44}) endows the form
$$
w_{k_1,j}(t)=m_1e^{ ta_{k_1,+}}-\myfrac{1}{a_{k\,+}}\myint{t}{0}e^{ (t-s)a_{k,+}}F_{k_1,j}(s)ds-\myint{-\infty}{t}F_{k_1,j}(s)ds.
$$
Because 
$$\myint{-\infty}{t}F_{k_1,j}(s)ds=O(e^t)\quad\text{as }\;t\to\infty,
$$
we conclude that for $k=k_1$, there holds
    $$
|w_{k_1,j}(t)|
\leq c_{27}\left\{\BA {lll}|t|e^t&\quad\text{ if }\, a_{k_1,+}=1,\\[1mm]
\max\{e^t,e^{ ta_{k_1,+}}\}&\quad\text{ if }\, a_{k_1,+}\neq 1
 \EA\right.$$
and finally  we infer (\ref{2-C43}), which complete the proof.\qeda\bigskip

\noindent {\bf Proof of  Theorem A.}  Assume that $u\in C^2(\overline\Gw\setminus\{0\})$ is a positive solution of $\CL_\gm u=0$ vanishing on 
 $\prt\Gw\setminus\{0\}$.\smallskip 

 {\it Case 1: $\gm>\gm_1$.}  We claim (\ref{2-C17}) holds for some $c_{13}\geq 0$. \smallskip
 
  By \rlemma{trace}, (\ref{2-C16}) holds for some $a>0$. If $a<-\ga_-$, then (\ref{2-C17}) holds with $c_{13}=0$. If $a=-\ga_-$, then 
(\ref{2-C17}) holds by \rprop{trace2}-(i). Hence we are left with the case $a>-\ga_-$. As in the proof of \rprop{trace2} we define 
$k_1$ and $k_2$. By replacing $a$ by $a'=a+\ge$, we can assume that $a_{k_2,+}\neq 1$ and $a_{k_1,-}\neq 1$, to avoid the resonance complication in (\ref{2-C39}) and (\ref{2-C43}), hence
   $$
\norm{v(t,.)}_{C^1(\overline{\BBS^{N-1}_+})}
\leq c_{27}\left(e^{ta_{k_2,+}}+e^{ ta_{k_1,-}}+e^t\right).
$$
Furthermore $k_2=k_1+1$ and 
$$a_{k_2\,+}-a_{k_1\,-}=\myfrac{1}{2}\left(\sqrt{4\gm+4\gl_{k_1+1}+(N-2)^2}+\sqrt{4\gm+4\gl_{k_1}+(N-2)^2}\right)>0,
$$
which yields
   $$
\norm{v(t,.)}_{C^1(\overline{\BBS^{N-1}_+})}
\leq c_{28}\left(e^{ ta_{k_1,-}}+e^t\right).
$$
This implies that $u$ satisfies  
   $$
u(x)
\leq c_{29}\left(|x|^{\ga_{k_1,-}}+|x|^{1-a}\right)\rho(x).
$$
We iterate this procedure up to obtain 
  $$
u(x)
\leq c_{30}|x|^{\ga_{-}}\rho(x)
$$
and we conclude as in the proof of \rprop{trace2}, Step 3.\qeda\medskip


{\it Case 2: $\gm=\gm_1$.} In this case, the difficulty comes from the fact that there is no dissipation of energy in 
(\ref{2-C26}) for $a=-\ga_-=\frac{N-2}{2}$. But from the above iterative procedure in the {\it Case 1}, we could obtain  could  obtain that for  some $\gd\in(0,\,1)$, 
$$
u(x)
\leq c_{31}|x|^{-\frac{N-2}{2}-\gd}\rho(x).
$$

 {\it 
We finally show that  there exists $c_{32}\geq 0$ such that 
      \bel{2-C47}\displaystyle
\lim_{r\to 0} r^{\frac{N-2}{2}}\myfrac{u(r,.)}{\ln r}=-c_{32}\psi_1(.)
\ee
in $C^1(\BBS^{N-1}_+)$ and 
      \bel{2-C47x}\displaystyle
\lim_{r\to 0} r^{\frac{N}{2}}\myfrac{u_r(r,.)}{\ln r}=\frac{(N-2)c_{32}}{2}\psi_1(.)
\ee
uniformly in $\BBS^{N-1}_+$.}\smallskip

 Note that (\ref{2-C24}) reduces   that   
     $$
     \BA {lll}\displaystyle
(1+\ge_1(t))v_{tt}+\ge_2(t)v_t+(N-1+\ge_3(t))v+\Gd'v\\[2mm]
\phantom{-----}
+\langle \nabla' v,\ge_4(t,.)\rangle+\langle \nabla' v_t,\ge_5(t,.)\rangle+\langle \nabla'(\langle\nabla' v,{\bf e}_N\rangle),\ge_6(t,.)\rangle
=0,
\EA$$
in $(-\infty,r_0)\ti \BBS^{N-1}_+$, vanishes on $(-\infty,r_0)\ti\prt \BBS^{N-1}_+$ and the $\ge_j$ verify (\ref{2-C25}), and
   $$
v(t,\gs)\leq c_{33}e^{-\gd t}.
$$
Since the operator involved in the equation is uniformly elliptic we have by standard regularity theory
   $$
   \BA {lll}\displaystyle
\norm{v}_{C^{2,\gd}([T-1,T+1]\ti\overline{\BBS^{N-1}_+})}+\norm{v_t}_{C^{1,\gd}([T-1,T+1]\ti\overline{\BBS^{N-1}_+})}
+\norm{v_{tt}}_{C^{\gd}([T-1,T+1]\ti\overline{\BBS^{N-1}_+})}\\[3mm]
\phantom{------------------}\leq c_{34}\norm v_{L^\infty((T-2,T+2)\ti \BBS^{N-1}_+} \\[2mm]
\phantom{------------------}\leq c_{35} e^{-\gd T}
\EA$$
for any $T\leq r_0+3$. We set
$$X(t)=\myint{\BBS^{N-1}_+}{}v(t,.)\psi_1 dS,
$$
then 
\bel{2-C51}X''(t)+F(t)=0\ee 
where 
$$F(t)=\myint{\BBS^{N-1}_+}{}\!\!\left(\ge_1v_{tt}+\ge_2v_{t}+\ge_3 v+\langle \nabla' v,\ge_4(t,.)\rangle+\langle \nabla' v_t,\ge_5(t,.)\rangle+\langle \nabla'(\langle\nabla' v,{\bf e}_N\rangle),\ge_6(t,.)\rangle\right)\psi_1 dS.
$$
Hence 
\bel{2-C52}|F(t)|\leq c_{36}e^{(1-\gd)t}.
\ee
 
This implies that $X'(t)$ admits a limit $c_{37}\leq 0$ when $t\to-\infty$ and 
   $$
\lim_{t\to-\infty}t^{-1}X(t)=c_{37}. 
$$
Set 
$$W_2=\underset{k\geq 2}{\oplus}\ker(\Gd'+\gl_kId),
$$
and denote by $v_2$ the orthogonal projection of $v$ onto $W_2$. Then
     \bel{2-C54}\BA {lll}\displaystyle
v_{2\,tt}+(N-1)v_2+\Gd'v_2=F_2(t,.),
\EA\ee
where 
$$|F_2(t,.)|\leq c_{38}e^{(1-\gd)t}.
$$
Since $\gl_2=2N$, the function $Y(t)=\norm{v_2(t,.)}_{L^2(\BBS^{N-1}_+)}$ satisfies in $\CD'(-\infty,r_1)$
  $$
Y''-(N+1)Y\geq -c_{38}e^{(1-\gd)t}.
$$
Because $Y(t)=o(e^{-\sqrt{N+1}t})$ when $t\to-\infty$, it follows by the maximum principle that $Y(t)=O(e^{\sqrt{N+1}t}+e^{(1-\gd)t})=O(e^{(1-\gd)t})$. Using again the standard regularity estimates for elliptic equations, we derive
     \bel{2-C56}\BA {lll}\displaystyle
\norm {v_2(t,.)}_{C^1(\overline{\BBS^{N-1}_+})}+\norm {v_{2\,t}(t,.)}_{C(\overline{\BBS^{N-1}_+})}\leq c_{39}e^{(1-\gd)t}.
\EA\ee
Combining (\ref{2-C51}) and (\ref{2-C54}) we derive (\ref{2-C47}). Since $v(t,.)=X(t)\psi_1+v_2(t,.)$ it follows from (\ref{2-C56}) that 
    $$
\lim_{t\to-\infty}v_t(t,.)=c_{37}\psi_1\quad \text{uniformly in } \BBS^{N-1}_+.
$$
  Thus, the indentity $u_r(r,\cdot)=r^{-\frac{N}{2}}\left(\frac{2-N}{2}v(t,\cdot)+v_t(t,\cdot)\right)$    implies (\ref{2-C47}) and (\ref{2-C47x}).  \hfill\qeda 

\subsection{ Existence and uniqueness }

 \noindent{\bf Proof of Theorem B. } We still assume that $\Gw$ satisfies the condition $(\CC\text{-}1)$ and $\prt\BBR^N_+$ is tangent to $\prt\Gw$ at $0$. For $\ge>0$ let $u_\ge$ be the solution of 
 \bel{2-D5}\left\{\BA {lll}
\CL_\gm u_\ge=0\quad&\text{in }\ \Gw_\ge:=\Gw\setminus \overline B_\ge,\\[0.5mm]
\phantom{\CL_\gm }u_\ge=0&\text{on }\, \prt\Gw\cap\overline B^c_\ge,\\[0.5mm]
\phantom{\CL_\gm }u_\ge=\phi_\gm&\text{on }\, \Gw\cap \prt B_\ge.
\EA\right.\ee
Since $\Gw\subset\BBR^N_+$, $u_\ge\leq \phi_\gm$ in $\Gw_\ge$ and 
 \bel{2-D6}\left.\myfrac{\prt u_\ge}{\prt{\bf n}}\right|_{\Gw\cap \prt B_\ge}\leq \left.\myfrac{\prt \phi_\gm}{\prt{\bf n}}\right|_{\Gw\cap \prt B_\ge}<0,
\ee
where $\bf n=\ge^{-1}x$. Furthermore, if $0<\ge'<\ge$, $u_{\ge'}\lfloor_{\Gw\cap \prt B_\ge}\leq u_{\ge}\lfloor_{\Gw\cap \prt B_\ge}=\phi_\gm$, hence $u_{\ge'}\leq u_{\ge}$ in $\Gw_\ge$. There exists $u_0=\lim_{\ge\to 0}u_\ge$ and $u_0$ is a nonnegative solution of $\CL_\gm u=0$ in $\Gw$ which vanishes on $\prt\Gw\setminus\{0\}$ and is smaller than $\phi_\gm$. 

Let 
$\gz\in\BBX_\gm(\Gw)$, $\gz>0$, then, with ${\bf n}'=-\frac x{|x|}=-{\bf n}$,
$$
\BA {lll}0=\myint{\Gw_\ge}{}\gz\gg^\Gw_\gm\CL_\gm u_\ge dx\\[4mm]\phantom{0}
=\myint{\Gw_\ge}{}u_\ge\CL^*_\gm \gz d\gg^\Gw_\gm
+\myint{\prt B_\ge\cap \Gw}{}\left(-\myfrac{\prt u_\ge}{\prt{\bf n}'}\gz\gg^\Gw_\gm+\left(\gz\myfrac{\prt \gg^\Gw_\gm}{\prt{\bf n}'}+\gg^\Gw_\gm\myfrac{\prt \gz}{\prt{\bf n}'}\right)u_\ge\right) dS.
\EA$$
Using (\ref{2-D5}) and (\ref{2-D6}) we obtain
 $$
 \BA {lll}\myint{\Gw_\ge}{}u_\ge\CL^*_\gm \gz d\gg^\Gw_\gm
\geq \myint{\prt B_\ge\cap \Gw}{}\left(\myfrac{\prt \phi_\gm^\Gw}{\prt{\bf n}}\gz\gg^\Gw_\gm-\left(\gz\myfrac{\prt \gg^\Gw_\gm}{\prt{\bf n}}+\gg^\Gw_\gm\myfrac{\prt \gz}{\prt{\bf n}}\right)\phi_\gm^\Gw\right) dS.
\EA$$
We take $\gz=1$, hence $\CL^*_\gm \gz=\ell_\gm^\Gw$ and we get
 $$
 \BA {lll}\ell_\gm^\Gw\myint{\Gw_\ge}{}u_\ge d\gg^\Gw_\gm
\geq \myint{\prt B_\ge\cap \Gw}{}\left(\myfrac{\prt \phi_\gm^\Gw}{\prt{\bf n}}\gg^\Gw_\gm-\myfrac{\prt \gg^\Gw_\gm}{\prt{\bf n}}\phi_\gm^\Gw\right) dS\\[4mm]\phantom{\ell_\gm^\Gw\myint{\Gw_\ge}{}u_\ge d\gg^\Gw_\gm}
\geq 2\sqrt{\gm+\gm_1}\myint{\BBS_+^{N-1}}{}\psi_1^2dS-o (1),
\EA$$
 in the case $\gm>\gm_1$, and
$$\BA {lll}\ell_\gm^\Gw\myint{\Gw_\ge}{}u_\ge d\gg^\Gw_\gm
\geq \myint{\prt B_\ge\cap \Gw}{}\left(\myfrac{\prt \phi_\gm^\Gw}{\prt{\bf n}}\gg^\Gw_\gm-\myfrac{\prt \gg^\Gw_\gm}{\prt{\bf n}}\phi_\gm^\Gw\right) dS\\[4mm]\phantom{\ell_\gm^\Gw\myint{\Gw_\ge}{}u_\ge d\gg^\Gw_\gm}
\geq \left(\myfrac N2-1\right)\myint{\BBS_+^{N-1}}{}\psi_1^2dS-o (1),
\EA$$
 in the case $\gm=\gm_1$. Since $u_\ge\leq \phi_\gm^\Gw$, 
 $$u_\ge \gg^\Gw_\gm\leq\gg^\Gw_\gm\phi_\gm^\Gw=r^{2-N}\psi_1^2\in L^1(\Gw).$$
 Therefore, by dominated convergence theorem, we conclude that
\bel{2-D10}\ell_\gm^\Gw\myint{\Gw}{}u_0 d\gg^\Gw_\gm\geq\left\{\BA {lll}2\sqrt{\gm+\gm_1}\myint{\BBS_+^{N-1}}{}\psi_1^2dS\quad&\text{if }\;
 \gm>\gm_1,\\[4mm]
  \left(\myfrac N2-1\right)\myint{\BBS_+^{N-1}}{}\psi_1^2dS\quad&\text{if }\;
 \gm=\gm_1.
 \EA\right.
\ee
 We infer that the function $u_0$ is nonzero. It is a positive solution of $\CL_\gm u_0=0$ in $\Gw$ which vanishes on $\prt\Gw\setminus\{0\}$. It follows from {\bf Theorem A} that there exists $k\geq 0$ such that 
 $$ \BA {lll}\displaystyle\lim_{x\to0} \frac{u(x)}{\rho(x)|x|^{\alpha_--1}}=k \quad&\text{if }\;
 \gm>\gm_1,\\[4mm]
 \displaystyle \lim_{x\to0} \frac{u(x)}{\rho(x)|x|^{-N/2}\ln |x|}=k \quad&\text{if }\;
 \gm=\gm_1.
 \EA 
 $$
Next we next show that $k=1$. In fact, if $k<1$, there exists $\ge_0>0$ such that for any $\ge\in(0,\ge_0)$
$$u_\ge\leq \frac{k+1}{2} \phi_\gm^\Gw,$$ 
 and then
 $$\lim_{\ge\to0^+}\ell_\gm^\Gw \myint{\Gw_\ge}{}u_\ge d\gg^\Gw_\gm\leq \frac{k+1}{2}
\ell_\gm^\Gw \myint{\Gw_\ge}{}\phi_\gm^\Gw d\gg^\Gw_\gm<\left\{\BA {lll}2\sqrt{\gm+\gm_1}\myint{\BBS_+^{N-1}}{}\psi_1^2dS\quad&\text{if }\;
 \gm>\gm_1,\\[4mm]
  \left(\myfrac N2-1\right)\myint{\BBS_+^{N-1}}{}\psi_1^2dS\quad&\text{if }\;
 \gm=\gm_1,
 \EA\right. $$ 
 which contradicts (\ref{2-D10}). Thus, (\ref{1-A25}) and (\ref{1-A26}) hold true. \qeda\medskip\medskip

 \noindent{\bf Proof of Corollary C.}
 {\it Identity (\ref{1-Ax32z}). } As a consequence of \rprop{trace2}, for any $\gz\in\BBX_\gm(\Gw)$ and $\ge>0$ we set $\Gw_\ge=\Gw\cap\overline B_\ge^c$, and  there holds
 $$\BA {lll}0=\myint{\Gw^\ge}{}\gz\gg^\Gw_\gm\CL_\gm\phi_\gm^\Gw dx\\[4mm]\phantom{0}
 =\myint{\Gw^\ge}{}\phi_\gm^\Gw\CL^*_\gm\gz d\gg^\Gw_\gm+
 \myint{\Gw\cap \prt B_\ge}{}\left(-\myfrac{\prt \phi_\gm^\Gw}{\prt{\bf n}}\gz\gg^\Gw_\gm+
 \left(\gg^\Gw_\gm\myfrac{\prt \gz}{\prt{\bf n}}+\gz\myfrac{\prt \gg^\Gw_\gm}{\prt{\bf n}}\right)\phi_\gm^\Gw\right)dS.
 \EA$$
 Using {\bf Proposition A.1} we have
 $$\displaystyle \myint{\Gw\cap \prt B_\ge}{}\left(-\myfrac{\prt \phi_\gm^\Gw}{\prt{\bf n}}\gz\gg^\Gw_\gm+
 \left(\gg^\Gw_\gm\myfrac{\prt \gz}{\prt{\bf n}}+\gz\myfrac{\prt \gg^\Gw_\gm}{\prt{\bf n}}\right)\phi_\gm^\Gw\right)dS
 =-\gz(0)A(\ge)(1+o (1)),
 $$
 where $A(\ge)$ is defined in (\ref{Ae}).
 
 The uniqueness follows direct from Kato's inequality (\ref{k2}). \qeda

 \mysection{The Dirichlet problem}
 
 \noindent{\bf Proof of Theorem D.} Note that in section \S3.2 for $\gl \in \frak M(\prt\Gw; \gb _\gm)$, problem 
$$
\left\{\BA {lll}
\CL_\gm u=0\quad&\text{in }\ \Gw,\\[0.5mm]
\phantom{\CL_\gm}
u=\gl\quad&\text{on }\,\prt\Gw
\EA\right.$$
 has a unique solution, denoting  $\BBK_\mu^\Omega(\lambda)$, which verifies the indentity
 $$\int_\Gw \BBK_\mu^\Omega(\lambda) \CL^*_\mu \gz d\gg^\Gw_\mu=\int_{\prt \Gw}\gz d(\gl\gb^\Gw_\mu)\quad\text{for all }\gz\in \BBX_\mu(\Gw).$$

Moreover, problem 
$$
\left\{\BA {lll}
\CL_\gm u=\nu\quad&\text{in }\ \Gw,\\[0.5mm]
\phantom{\CL_\gm}
u=0\quad&\text{on }\,\prt\Gw
\EA\right.$$
 has a unique solution, denoting  $\BBG_\mu^\Omega(\nu)$, which verifies the indentity
 $$\int_\Gw \BBK_\mu^\Omega(\lambda) \CL^*_\mu \gz d\gg^\Gw_\mu=\int_{ \Gw}\gz d\gg^\Gw_\mu\quad\text{for all }\gz\in \BBX_\mu(\Gw).$$ 
  Together with {\bf Corollary C} and the linearity of operator $\CL_\mu$, we have that $\BBK_\mu^\Omega(\lambda)+\BBG_\mu^\Omega(\nu)+k\phi^\mu_\Gw$ is a weak solution of (\ref{1-B1}) satisfying (\ref{1-Ax32}) and the uniqueness follows directly from Kato's inequality (\ref{k2}). \qeda\bigskip

 Our final part is to classify the boundary data for nonnegative $\CL_\gm$-harmonic function.\medskip
 
\noindent{\bf Proof of Theorem E. }    Let $\Gw$ be a bounded $C^2$ domain and $u$ be a nonnegative $\CL_\gm$-harmonic function in $\Gw$. 
We now show that there exists a nonnegative measure $\gl$ on $\prt\Gw\setminus\{0\}$ and $k\geq 0$ such that 
   \bel{2-D-12}
u=\BBK^\Gw_\gm[\gl]+k\phi_\gm^\Gw.
  \ee
 
 For $\ge>0$ the term $\gm|x|^{-2}$ is bounded in $\Gw_\ge=\Gw\cap\overline B_\ge^c$. Hence the exists a nonnegative Radon measure $\gl_\ge$ such that $u$ is the unique solution of 
   $$
   \left\{\BA {lll}
\CL_\gm u=0\quad&\text{in }\ \Gw_\ge,\\[0.8mm]
\phantom{\CL_\gm}u=\gl_\ge\quad&\text{on }\;\prt\Gw_\ge.
\EA \right. $$
Furthermore $\gl_\ge$ is the boundary trace is achieved in {\it dynamical sense}, see \cite{MaVe0} and references therein. Hence for any 
$\gz\in C(\overline\Gw)$ vanishing on $ B_\ge$, there holds
   $$
\lim_{\gd\to 0}\myint{\Gs_\gd}{}u\gz dS=\myint{\prt \Gw\cap B_\ge^c}{}\gz d\gl_\ge,
$$
where $\Gs_\gd=\{x\in\Gw:\gr(x)=\gd\}$. If we write 
$$\gl_\ge=\gl_\ge\lfloor_{\prt \Gw\cap B_\ge^c}+u\lfloor_{\Gw\cap \prt B_\ge},$$
it proves that for $0<\ge'<\ge$, $\gl_\ge\lfloor_{\prt \Gw\cap B_\ge^c}=\gl_{\ge'}\lfloor_{\prt \Gw\cap B_\ge^c}$. This defines in a unique way a nonnegative Radon $\gl$ on $\prt\Gw\setminus\{0\}$ measure such that (\ref{2-D-12}) holds for all $\gz\in \BBX_\gm(\Gw)$ vanishing near $0$. Furthermore $\gr u\in L^1(\Gw_\ge)$ for any $\ge>0$. Denote by $\BBK^{\Gw_\ge}_\gm$ the Poisson potential of $\CL_\gm$ in $\Gw_\ge$. Then
$$
 u\lfloor_{\Gw_\ge}=\BBK^{\Gw_\ge}_\gm[\gl_\ge\lfloor_{\prt \Gw\cap B_\ge^c}]+\BBK^{\Gw_\ge}_\gm[u\lfloor_{\Gw\cap\prt B_\ge}].
$$
For $0<\ge'<\ge$, one has that $\BBK^{\Gw_{\ge'}}_\gm[\gl_{\ge'}\lfloor_{\prt \Gw\cap B_{\ge'}^c}]\lfloor_{\Gw\cap\prt B_\ge}\geq 0$. Therefore
$\BBK^{\Gw_{\ge'}}_\gm[\gl_{\ge'}\lfloor_{\prt \Gw\cap B_{\ge'}^c}]\geq \BBK^{\Gw_\ge}_\gm[\gl_\ge\lfloor_{\prt \Gw\cap B_\ge^c}]$ in $\Gw_\ge$. Hence 
  $$
\lim_{\ge\to 0}\BBK^{\Gw_\ge}_\gm[\gl_\ge\lfloor_{\prt \Gw\cap B_\ge^c}]=\BBK^{\Gw}_\gm[\gl]\leq u\quad\text{in }\;\Gw.
$$

Next we aim to characterize the behaviour at $0$. By contradiction we assume that 
    $$
\limsup_{\ge\to 0}\myint{\Gw\cap \prt B_\ge}{}ud\gb_\gm^\Gw =\lim_{\ge_k\to 0}\myint{\Gw\cap \prt B_{\ge_k}}{}ud\gb_\gm^\Gw=\infty. 
$$
Then for any $m>0$ there exists a sequence $\{\ge_{m,k}\}\subset\BBR_+^*$ tending to $0$ and a sequence $\{\ell_{m,k}\}\subset\BBR_+^*$ tending to $\infty $ such that 
$$
\myint{\Gw\cap \prt B_{\ge_k}}{}\min\{u,\ell_{m,k}\}d\gb_\gm^\Gw =m.
$$
Set $\gt_{m,k}=\min\{u\lfloor_{\Gw\cap \prt B_{\ge_k}},\ell_{m,k}\}$ and set $u_{m,k}=\BBK^{\Gw_{\ge_k}}_\gm[\gt_{m,k}\chi_{_{\Gw\cap \prt B_{\ge_k}}}]$.
Then
$$u_{m,k}\leq u\qquad\text{in }\;\Gw_{\ge_k},
$$
and we recall that
$$\myint{\Gw\cap \prt B_{\ge_k}}{}\phi_\gm^\Gw d\gb_\gm^\Gw =c_\gm (1+\circ (1)),
$$
where $c_\gm$ is the constant defined in (\ref{2-A7}).  Combining the boundary Harnack inequality with the standard Harnack inequality, one infers
 \bel{2-D-17}
 c_{47}\myfrac{\gr(x)}{\gr(y)}\leq
c_{46}\myfrac{\phi_\gm^\Gw(x)}{\phi_\gm^\Gw(y)}\leq \myfrac{u_{m,k}(x)}{u_{m,k}(y)}\leq c_{44}\myfrac{\phi_\gm^\Gw(x)}{\phi_\gm^\Gw(y)}\leq c_{45}\myfrac{\gr(x)}{\gr(y)}
\ee
for all $x,y\in\Gw$ such that $|x|=|y|\geq 2\ge_k$. If we set 
$\dot\phi_\gm^\Gw(x)=\myfrac{\phi_\gm^\Gw(x)}{\gr(x)}$ and $\dot u_{m,k}(x)=\myfrac{u_{m,k}(x)}{\gr(x)}$,
then (\ref{2-D-17}) becomes
  \bel{2-D-17'}
 c_{47}\leq
c_{46}\myfrac{\dot\phi_\gm^\Gw(x)}{\dot\phi_\gm^\Gw(y)}\leq \myfrac{\dot u_{m,k}(x)}{\dot u_{m,k}(y)}\leq c_{44}\myfrac{\dot\phi_\gm^\Gw(x)}{\dot\phi_\gm^\Gw(y)}\leq c_{45}.
 \ee

Assume for a while that we have proved that there exists $\gth>0$, independent of $m$ and $k$ such that for for any 
     \bel{2-D-18}\BA {lll}\displaystyle
\myint{\prt B_{2\ge_k}\cap\Gw}{}\dot u_{m,k}d\left(\gr\gb_\gm^\Gw\right)\geq\gth \myint{\prt B_{\ge_k}\cap\Gw}{}\dot u_{m,k}d\left(\gr\gb_\gm^\Gw\right)=\gth m.
\EA  \ee
If we assume that for $\gd\leq 2\ge_{k_0}$
$$2c_\gm\geq\myint{\Gw\cap \prt B_{\gd}}{}\dot \phi_\gm^\Gw d\left(\gr\gb_\gm^\Gw\right) \geq\frac{c_\gm}{2},
$$
one has for $k\geq k_0$,
$$\myint{\Gw\cap \prt B_{2\ge_k}}{}\dot u_{m,k} d\left(\gr\gb_\gm^\Gw\right) \geq \gth m\geq \myfrac{\gth m}{2c_\gm}\myint{\Gw\cap \prt B_{2\ge_k}}{}\dot \phi_\gm^\Gw d\left(\gr\gb_\gm^\Gw\right)
$$
Since
$$\dot \phi_\gm^\Gw (x)\leq \myfrac{c_{45}}{c_{46}}\dot \phi_\gm^\Gw (y)
$$
and
$$\dot u_{m,k}(x)\geq c_{47}\dot u_{m,k}(y),
$$
we derive
  $$
\myfrac{1}{c_{47}}\dot u_{m,k}(x)\myint{\prt B_{2\ge_k}\cap\Gw}{}d\left(\gr\gb_\gm^\Gw\right)\geq\myfrac{\gth mc_{46}}{2c_\gm c_{45}}\dot \phi_\gm^\Gw (x)\myint{\Gw\cap \prt B_{2\ge_k}}{}d\left(\gr\gb_\gm^\Gw\right).
$$
Therefore 
$$u_{m,k}(x)\geq c_{48}m\phi_\gm^\Gw (x)\quad\text{for all }\;x\in\Gw\,\text{ s.t. }\abs x=2\ge_k,
$$
and $c_{48}>0$ is independent of $m$ and $\ge_k$. This implies by the maximum principle and letting $\ge_k\to 0$
    \bel{2-D-20}\BA {lll}\displaystyle
u(x)\geq u_{m,k}(x)\geq c_{48}m\phi_\gm^\Gw (x)\quad\text{for all }\;x\in\Gw. \EA  \ee
Since $m$ is arbitrary we obtain a contradiction. Hence there holds
\bel{2-D-21}\displaystyle
\limsup_{\ge\to 0}\myint{\Gw\cap \prt B_\ge}{}ud\gb_\gm^\Gw =\lim_{\ge_k\to 0}\myint{\Gw\cap \prt B_{\ge_k}}{}ud\gb_\gm^\Gw=m_u<\infty. 
\ee
Then inequality (\ref{2-D-20}) holds without truncation with $m$ replaced by $m_u$. We recall that
\bel{2-D-21x}
w_\ge:=\BBK^{\Gw_\ge}_\gm[u\lfloor_{\Gw\cap\prt B_\ge}]=u\lfloor_{\Gw_\ge}-\BBK^{\Gw_\ge}_\gm[\gl_\ge\lfloor_{\prt \Gw\cap B_\ge^c}]\quad\text{in }\,\Gw_\ge.
\ee
{\it Case 1: We first assume that $m_u>0$}. Then (\ref{2-D-18}) combined with the maximum principle yields
   $$
\myint{\prt B_{\ge_k}\cap\Gw}{}\dot w_{\ge_k}d\left(\gr\gb_\gm^\Gw\right)\geq \myint{\prt B_{2\ge_k}\cap\Gw}{}\dot w_{\ge_k}d\left(\gr\gb_\gm^\Gw\right)\geq\gth \myint{\prt B_{\ge_k}\cap\Gw}{}\dot w_{\ge_k}d\left(\gr\gb_\gm^\Gw\right)=\gth m_u(1+o(1))
$$
with $\dot w_\ge=\gr^{-1}w_\ge$.  Inequality (\ref{2-D-17'}) is replaced by 
\bel{2-D-22}
 c_{47}\leq
c_{46}\myfrac{\dot\phi_\gm^\Gw(x)}{\dot\phi_\gm^\Gw(y)}\leq \myfrac{\dot w_\ge(x)}{\dot w_\ge(y)}\leq c_{44}\myfrac{\dot\phi_\gm^\Gw(x)}{\dot\phi_\gm^\Gw(y)}\leq c_{45}\quad\text{in }\,\Gw_{2\ge}.
\ee
Therefore, for $\ge_k$ small enough and $|x|=2\ge_k$,
$$\BA {lll}
\dot w_\ge(x)\myint{\prt B_{\ge_k}\cap\Gw}{}d\left(\gr\gb_\gm^\Gw\right)\leq c_{45}\myint{\prt B_{\ge_k}\cap\Gw}{}\dot w_{\ge_k}(y)d\left(\gr\gb_\gm^\Gw\right)\leq 2c_{45}m_u\\[4mm]
\phantom{\dot w_\ge(x)\myint{\prt B_{\ge_k}\cap\Gw}{}d\left(\gr\gb_\gm^\Gw\right)--}\leq 
\myfrac{4c_{45}m_u}{c_\gm}\myint{\prt B_{\ge_k}\cap\Gw}{}\dot \phi_\gm^\Gw(y)d\left(\gr\gb_\gm^\Gw\right)
\leq \myfrac{4c_{44}c_{45}m_u}{c_\gm c_{47}}\dot \phi_\gm^\Gw(x)\myint{\prt B_{\ge_k}\cap\Gw}{}d\left(\gr\gb_\gm^\Gw\right),
\EA$$
which implies 
\bel{2-D-23x}
w_{\ge_k}(x)\leq \myfrac{4c_{44}c_{45}m_u}{c_\gm c_{47}}\phi_\gm^\Gw(x):=c_{49}m_u\phi_\gm^\Gw(x)\quad\text{for }\;x\in\Gw\cap\prt B_{2\ge_k}.
 \ee
 Hence
 \bel{2-D-24}
w_{\ge_k}(x)\leq \myfrac{4c_{44}c_{45}m_u}{c_\gm c_{47}}\phi_\gm^\Gw(x):=c_{49}m_u\phi_\gm^\Gw(x)\quad\text{for }\;x\in\Gw\cap\prt B_{2\ge_k}.
 \ee
 Since $w_{\ge_k}$ and $\phi_\gm$ are $\CL_\gm$-harmonic in $\Gw_{2\ge_k}$, and vanishes on $\prt\Gw\cap B_{2\ge_k}$ it follows that inequality 
 (\ref{2-D-24}) also holds for any $x\in \Gw_{2\ge_k}$. By definition $w_{\ge_k}=\BBK_\gm^{\Gw_{\ge_k}}[u\lfloor_{\Gw\cap\prt B_{\ge_k}}]$, hence 
  \bel{2-D-25}
\BBK_\gm^{\Gw_{\ge_k}}[u\lfloor_{\Gw\cap\prt B_{\ge_k}}](x)\leq c_{49}m_u\phi_\gm^\Gw(x)\quad\text{for }\;x\in\Gw_{2\ge_k}.
 \ee
Next we obtain the estimate from below. From (\ref{2-D-22}), with $|x|=2\ge_k$,
$$\BA {lll}
\dot w_{\ge_k}(x)\myint{\prt B_{\ge_k}\cap\Gw}{}d\left(\gr\gb_\gm^\Gw\right)\geq \ c_{47}\myint{\prt B_{\ge_k}\cap\Gw}{}\dot w_{\ge_k}(y)d\left(\gr\gb_\gm^\Gw\right) \ \geq \myfrac{c_{47}m_u}{2}\\[4mm]
\phantom{\dot w_\ge(x)\myint{\prt B_{\ge_k}\cap\Gw}{}d\left(\gr\gb_\gm^\Gw\right)\ }\geq 
\myfrac{c_{47}m_u}{4c_\gm}\myint{\prt B_{\ge_k}\cap\Gw}{}\dot \phi_\gm^\Gw(y)d\left(\gr\gb_\gm^\Gw\right)
\geq \myfrac{c_{47}c_{44}m_u}{4c_\gm c_{45}}\dot \phi_\gm^\Gw(x)\myint{\prt B_{\ge_k}\cap\Gw}{}d\left(\gr\gb_\gm^\Gw\right).
\EA$$
Hence 
$$
w_{\ge_k}(x)\geq \myfrac{c_{47}c_{44}m_u}{4c_\gm c_{45}}\phi_\gm^\Gw(x):=c_{50}m_u\phi_\gm^\Gw(x)\quad\text{for }\;x\in\Gw\cap\prt B_{2\ge_k}.
$$
 It follows that
   \bel{2-D-27}
\BBK_\gm^{\Gw_{\ge_k}}[u\lfloor_{\Gw\cap\prt B_{\ge_k}}](x)\geq c_{50}m_u\phi_\gm^\Gw(x)\quad\text{for }\;x\in\Gw_{2\ge_k}.
 \ee
 From (\ref{2-D-21x}), (\ref{2-D-27}) and (\ref{2-D-25}) we infer
 \bel{2-D-28}
c_{50}m_u\phi_\gm^\Gw\leq u\lfloor_{\Gw_{\ge_k}}-\BBK^{\Gw_{\ge_k}}_\gm[\gl_{\ge_k}\lfloor_{\prt \Gw\cap B_{\ge_k}^c}]\leq c_{48}m_u\phi_\gm^\Gw\quad\text{in }\,\Gw_{2\ge_k}.
\ee
This implies, by letting $\ge_k\to 0$, 
$$
c_{50}m_u\phi_\gm^\Gw\leq u-\BBK^{\Gw}_\gm[\gl]\leq c_{48}m_u\phi_\gm^\Gw\quad\text{in }\,\Gw.
$$
Therefore, the function $u-\BBK^{\Gw}_\gm[\gl]$ is $\CL_\gm$-harmonic and positive in $\Gw$ and it vanishes on $\prt\Gw$. By Corollary C, it implies that it coincides with $c\phi_\gm^\Gw$ for some $c\geq 0$ (and in that case $c_{50}m_u\leq c\leq c_{49}m_u$). \smallskip

\nind {\it Case 2: Assume $m_u=0$}. Following the same inequalities as in Case 1, (\ref{2-D-23x}) is replaced by: for any $\gd>0$ there exists $k_0>0$ such that for $k\geq k_0$, 
$$
w_{\ge_k}(x)\leq \gd\phi_\gm^\Gw(x)\quad\text{for }\;x\in\Gw\cap\prt B_{2\ge_k}.
 $$
 Hence (\ref{2-D-28}) is transformed into 
 $$
0\leq u\lfloor_{\Gw_{\ge_k}}-\BBK^{\Gw_{\ge_k}}_\gm[\gl_{\ge_k}\lfloor_{\prt \Gw\cap B_{\ge_k}^c}]\leq \gd\phi_\gm^\Gw\quad\text{in }\,\Gw_{2\ge_k}.
$$
Letting successively $\ge_k\to 0$ and $\gd\to 0$ yields $u-\BBK^{\Gw}_\gm[\gl]=0$ in $\Gw$, which ends the proof. \qeda

\section*{Appendix: Estimates (\ref{1-C1}) and (\ref{1-A23})}

 \noindent{\bf Proposition A.1} {\it Assume $\Gw$ is a bounded $C^2$ domain such that $0\in\prt \Gw$ satisfying condition (C-1) and let $\gg_\gm^\Gw$ be defined by (\ref{1-A22}) and normalized by $\norm{\gg_\gm^\Gw}_{L^2(\Gw)}=1$. Then 
     $$
\lim_{r\to 0}r^{1-\ga_+}\gg_\gm^\Gw(r,.)=c_1\psi_1\quad\text{in }\,C_{loc}^{1}(\BBS^{N-1}_+)
$$
and
$$
\lim_{r\to 0}r^{-\ga_+}\gg_{\gm\,r}^\Gw(r,.)=c_1\left(1-\myfrac{N}{2}\right)\psi_1\quad\text{locally uniformly in }\BBS^{N-1}_+.
$$
}
     
     
\nind {\bf Proof. }     Since  $\ga_++(N-2)\ga_+-\gm+1-N=0$, the function $x\mapsto w(x):=|x|^{\ga_+}$ satisfies
     $$\tilde\CL_\gm w(x):=\CL_\gm w-\ell_\gm^\Gw w=|x|^{\ga_+-2}\left(N-1-\ell_\gm^\Gw|x|^2\right)\quad\text{in }\,\BBR^N\setminus\{0\}.
     $$
 Furthermore,  $\nabla w\in L^2_{loc}(\BBR^N)$.
Let $R_0>0$ such that  $N-1\geq \ell_\gm^\Gw R_0^2$ and $m>0$ such that $mw\geq \gg_\gm^\Gw$ on $\Gw\cap B^c_{R_0}$. Then the function
$(\gg_\gm^\Gw-mw)_+$ belongs to $H_\mu(\Gw)$ and satisfies $\tilde\CL_\gm (\gg_\gm^\Gw-mw)_+\leq 0$ in the dual of $H_\mu(\Gw)$. Hence
$$\myint{\Gw}{}\left(|\nabla (\gg_\gm^\Gw-mw)_+|^2+\left(\myfrac{\gm}{|x^2|}-\ell_\gm^\Gw\right)(\gg_\gm^\Gw-mw)_+^2\right) dx\leq 0.
$$
Therefore $(\gg_\gm^\Gw-mw)_+\leq 0$, which implies that 
 $$
0< \gg_\gm^\Gw(x)\leq m|x|^{\ga_+}\quad\text{for all }\; x\in\Gw.
$$
Then we proceed as in the proof of \rprop{trace2}. We flatten the boundary near $0$ and set
$$ v(t,\gs)=r^{-\ga_+}\tilde \gg_\gm^\Gw(r,\gs)\quad{\rm with}\ \, t=\ln r,
$$ 
where the function $\tilde \gg_\gm^\Gw$ is defined similarly as  $\tilde u$ in (\ref{2-C19}). Then $v$ is bounded in $(-\infty,T_0]\ti \BBS^{N-1}_+$ where it satisfies 
$$
\BA {lll}\displaystyle
(1+\ge_1(t,.))v_{tt}+\left(N-2+2\ga_++\ge_2(t,.)\right) v_t+\left(\ga_+(\ga_++N-2)-\gm+\ge_3(t,.)+e^{2t}\ell_\gm^\Gw\right)v\\[2mm]
\phantom{---+e^{2t}\ell_\gm^\Gw}
+\Gd' v+\langle \nabla' v,\ge_4(t,.)\rangle+\langle \nabla' v_t,\ge_5(t,.)\rangle+\langle \nabla'(\langle\nabla' v,{\bf e}_N\rangle),\ge_6(t,.)\rangle
=0,
\EA$$
instead of (\ref{2-C24}). It vanishes on $(-\infty,T_0]\ti \prt \BBS^{N-1}_+$ and the $\ge_j$ satisfy again (\ref{2-C25}). \smallskip

{\it  Case 1: $\gm>\gm_1$.}  The energy method used in proof of \rprop{trace2} applies with no modification and we infer that there exists $c_{51}\geq 0$ such that 
$$
v(t,.)\to c_{51}\psi_1\quad\text{as }\; t\to -\infty
$$
in $C^1(\BBS^{N-1}_+)$ and $v_t(t,.)\to 0$ uniformly in $\BBS^{N-1}_+$.  If $c_{51}= 0$, we can prove, as in \rprop{trace2}-(ii) that there exists $\gt>0$ such that 
   \bel{3-A-4}
 \gg_\gm^\Gw(x)\leq c_{52}|x|^{\ga_++\gt}\quad\text{for all }\; x\in\Gw.
 \ee
Iterating this process, we infer that (\ref{3-A-4})  holds for any $\gt>0$. For $k>1$, let $\ga_{k,+}$ be the positive root of (\ref{2-A3}) and put 
$w_k(x)=|x|^{\ga_{k,+}}$. Then 
$$\tilde\CL_\gm w_k (x)=|x|^{\ga_{k,+}-2}\left(\gl_k-\ell_\gm^\Gw|x|^2\right)\quad\text{in }\,\BBR^N\setminus\{0\}.
$$
Since $\gl_k\to\infty$, as $k\to\infty$, we choose $k$ such that $\gl_k>\ell_\gm^\Gw ($diam$(\Gw))^2$. Hence $w_k$ is a supersolution of $\tilde\CL_\gm$. Because $ \gg_\gm^\Gw(x)=o (w_k(x))$ near $x=0$, it follows that $ \gg_\gm^\Gw(x)\leq \ge o(w_k(x))$ in $\Gw$ for any $x\in\Gw$. Hence $\gg_\gm^\Gw=0$, which is a contradiction. Finally it implies that $c_{51}> 0$, which yields (\ref{1-C1})-(i). Because the convergence holds in $C^1(S^{N-1}_+)$ and $v_t(t,.)\to 0$, we infer
 $$
\lim_{r\to 0} r^{1-\ga_+}\nabla \tilde\gg_\gm^\Gw(r,.)=c_{51}\left(\ga_+\psi_1{\bf e}+\nabla'\psi_1\right)
$$
where ${\bf e}=\frac x{|x|}$. This implies the claim. \smallskip

{\it  Case 2: $\gm=\gm_1$.} Set $v(t,.)=r^{\frac N2-1}u(r,.)$ with $t=\ln r$ and $X(t)=\myint{\BBS^{N-1}_+}{}v(t,.)\psi_1 dS$ and obtain again (\ref{2-C51}), where 
$F(t,.)$ satisfies (\ref{2-C52}).  Since  $X't)\to 0$ and $X$ is bounded, it follows that $X(t)$ admits a limit $c_{52}\geq 0$ when $t\to-\infty$. As in the proof of 
{\bf Theorem A}, we infer that 
 $$
\lim_{t\to -\infty} v(t,.)=c_{52}\psi_1\quad\text{in }\; C^1(\BBS^{N-1}_+)\quad\text{and }\;\lim_{t\to -\infty} v_t(t,.)=0\quad\text{uniformly in }\; \BBS^{N-1}_+.
$$
If $c_{52}= 0$ we derive a contradiction as in the first case.\qeda\medskip

\noindent{\bf Proof of (\ref{1-A23}).} Since $\rho^*\leq \frac1{l^\Gw_\gm}$, then comparison
principle implies that $\sigma^\Gw_\gm\geq \gg^\Gw_\gm$ in $\Gw$. Next we show 
$\gs^\Gw_\gm\leq c_2\gg^\Gw_\gm$ in $\Gw$. In fact, we only have to show this inequality 
holds in a neighborhood of the origin.

{\it Case 1: the boundary is flat at the origin}.  We first prove above inequality when $\Omega$ is 
flat in a neighborhood of the origin, i.e. $B'_R\ti[0,R)\subset   \Omega\subset \R^N_+$ for some $R>0$. 

For $\tau\in\R$, denote 
$$w_\tau(x)=|x|^\tau x_N\ \ \text{and } \ \tilde w_\tau(x)=|x|^\tau x_N^2\ \ \text{in } \, \R^N_+,$$
and direct calculation shows that 
$$\CL_\mu w_\tau(x)=[\mu-\tau(\tau+N)]|x|^{\tau-2}x_N,\quad \CL_\mu \tilde w_\tau(x)=[\mu-\tau(\tau+N+2)]|x|^{\tau-2}x_N^2-2|x|^{\tau}\ \ \text{for }\, x\in\R^N_+.$$
Let 
$$\overline u(x)=\left\{\BA {lll}
w_{\alpha_+-1}-\frac12\tilde w_{\alpha_+-1} \quad&\text{if }\ \alpha_+\geq 0,\\[0.5mm]
\phantom{ }
w_{\alpha_+-1}-\frac{N-2}{2(N+2)}w_{\alpha_+}-\frac{2}{N+2}\tilde w_{\alpha_+-1} \quad&\text{if }\, \ \alpha_+< 0,
\EA\right.
$$
by resetting $R\in(0,1]$ such that $\overline u>0$ in $\R^{N-1}\ti (0,R]$.

When $\alpha_+\geq 0$, we have that $\mu-(\alpha_+-1)(\alpha_+-1+N)=0$ and
$$\CL_\mu \overline u(x)=  \alpha_+|x|^{\alpha_+-3} x_N^2+ |x|^{\alpha_+-1}\geq  |x|^{\alpha_+-1},$$
thus
there exists $t_1>0$ such that $t_1\overline u\geq \gs^\Gw_\gm $ on $\Omega \cap (\R^{N-1}\ti \{R\})$ and
$$\CL_\mu (t_1\overline u)(x)\geq t_1 |x|^{\alpha_+-1}\geq \CL_\mu \gs^\Gw_\gm(x)\quad\text{in }\,\Omega \cap (\R^{N-1}\ti (0,R))$$
By comparison principle, we have that 
$$\gs^\Gw_\gm\leq t_1\overline u$$
which, together with the inequality $\overline u\leq  2t_1w_{\alpha_+-1}$, implies that $\gs^\Gw_\gm\leq c_2\gg^\Gw_\gm$. 

When $\alpha_+\in[\frac{2-N}{2}, 0)$ if $N\geq 3$,
$$\BA {lll}
\CL_\mu \overline u(x) =\frac{N-2}{2(N+2)} (2\alpha_++2+N)|x|^{\alpha_+-2} x_N
+\frac{4}{N+2} \alpha_+|x|^{\alpha_+-3} x_N^2+\frac{4}{N+2} |x|^{\alpha_+-1} \\[4mm]
\phantom{\CL_\mu \overline u(x) }\geq \left(\frac{N-2}{2(N+2)}(2\alpha_++2+N) -\frac{4\alpha_+}{N+2}\right)|x|^{\alpha_+-2} x_N+\frac{4}{N+2} |x|^{\alpha_+-1} 
\\[4mm]
\phantom{\CL_\mu \overline u(x) }\geq \frac{4}{N+2} |x|^{\alpha_+-1}.
\EA$$
The remaining of the proof is similar to the previous one and we omit it.\smallskip

\smallskip

{\it Case 2: the boundary is not flat at origin. }

We define the function $\Gth=(\Gth_1,...,\Gth_N)$ on $D_R$ by $y_j=\Gth_j(x)=x_j$ if $1\leq j\leq N-1$ 
and $y_N=\Gth_N(x)=x_N-\gth(x')$. Since $D\Gth(0)=Id$ we can assume that $\Gth$ is a diffeomorphism from $D_R$ onto $\Gth(D_R)$. 
We set
$$ 
u_1(x)=\overline u(y)\qquad\text{for all }y\in D^+_R=B'_R\ti[0,R).
$$
Then by (\ref{2-C20}) and (\ref{2-C21}), we have that  
$$(-\Delta u_1(x)+\frac{\mu}{|x|^2}) u_1(x)=(-\Delta \overline u(y)+\frac{\mu}{|y|^2} \overline u(y))+
O(|y|) \left(|D^2 \overline u(y) |+ \frac{\mu}{|y|^2}  \overline u(y)\right) $$
Then by resetting $R>0$ small and the calculation in {\it Case 1}, we have that 
$$\CL_\mu u_1(x)\geq c_{53} |x|^{\alpha_+-1},\quad\forall\, x\in \Gth^{-1}(D^+_R).$$
By Hopf's Lemma, there exists $t_2>0$ such that $t_2  u_1\geq \gs^\Gw_\gm $ on $\Gth^{-1}(\prt B'_R\ti[0,R)) $ and by compactness of  $\Gth^{-1}(\overline B'_R\ti \{R\})$, there exists $t_3>0$ such that $t_3  u_1\geq \gs^\Gw_\gm $ on $\Gth^{-1}(\overline B'_R\ti \{R\})$. 
Applying comparison principle, for some $t_4\geq \max\{t_2,t_3\}$, we have that $$\gs^\Gw_\gm\leq t_4\overline u\quad {\rm in}\ \Gth^{-1}(D^+_R)$$
and
 we have $\gs^\Gw_\gm\leq c_2\gg^\Gw_\gm$ near the orgin. \qeda\medskip

 \noindent{\bf Proposition A.2} {\it Under the assumption of Proposition A.1 there exists $c_{53}>0$ such that 
\bel{5.15}
\lim_{r\to 0}r^{1-\ga_+}\gs_\gm^\Gw(r,.)=c_{53}\psi_1\quad\text{in }\,C_{loc}^{1}(\BBS^{N-1}_+)
\ee
and
\bel{5.16}
\lim_{r\to 0}r^{-\ga_+}\gs_{\gm\,r}^\Gw(r,.)=c_{53}\left(1-\myfrac{N}{2}\right)\psi_1\quad\text{locally uniformly in }\,\BBS^{N-1}_+.
\ee
}\medskip

\nind {\bf Proof. }  We follow the proof of Proposition A.1, flattening the boundary near $0$ and defining a new function $\tilde \gs_{\gm}^\Gw$ as previously.  By (\ref{1-A23}) the function
$$v(t,\gs)=r^{-\ga_+}\tilde \gs_{\gm}^\Gw(r,\gs)\qquad\text{with }\; t=\ln r,
$$
is bounded and it satisfies 
     $$
\BA {lll}\displaystyle
(1+\ge_1(t,.))v_{tt}+\left(N-2+2\ga_++\ge_2(t,.)\right) v_t+\left(\ga_+(\ga_++N-2)-\gm+\ge_3(t,.)\right)v\\[2mm]
\phantom{---+e^{2t}\ell_\gm^\Gw}
+\Gd' v+\langle \nabla' v,\ge_4(t,.)\rangle+\langle \nabla' v_t,\ge_5(t,.)\rangle+\langle \nabla'(\langle\nabla' v,{\bf e}_N\rangle),\ge_6(t,.)\rangle
=e^tm(t,.),
\EA$$
instead of (\ref{2-C24}), where the function $m$ is bounded as well as its gradient. Then $v$ satisfies the same bounds as the ones in the proof of \rprop {trace2}. The only difference is that the energy estimate $(\ref{2-C26})$ is replaced by
\bel{5.17}\BA {lll}\displaystyle
\myint{\BBS^{N-1}_+}{} \left(N-2-2\ga_++  \ge_2-\myfrac{1}{2}\prt_t\ge_{1}\right)v_t^2 dS-\myfrac{1}{2}\myint{\BBS^{N-1}_+}{}\prt_t\ge_{3} v^2 dS\\[4mm]\phantom{--}
=\myfrac{d}{dt}\left[\myint{\BBS^{N-1}_+}{}\left(\myfrac{1}{2}|\nabla v|^2-\myfrac{1}{2}\left[\ga_+(\ga_++2-N)-\mu+\ge_3\right]v^2-\myfrac{1}{2}(1+\ge_1)v_t^2 \right) dS\right]\\[4mm]
\phantom{--} -\myint{\BBS^{N-1}_+}{}\left(\langle\nabla' v,\ge_4\rangle+\langle\nabla' v_t,\ge_5\rangle+
\langle\nabla'(\langle\nabla' v,{\bf e}_N\rangle),\ge_6\rangle\right)v_t^2 dS+e^t\myint{\BBS^{N-1}_+}{}m(t,.)v(t,.)dS.
 \EA\ee
Therefore, if $2\ga_+\neq N-2$, we conclude that $(\ref{2-C28})$ holds, and $(\ref{2-C29})$ follows. We infer $(\ref{5.15})$ and $(\ref{5.16})$ as in the proof of \rprop {trace2}. When $2\ga_+= N-2$ the proof of $(\ref{2-C28})$ and $(\ref{2-C28})$ in the case $2\ga_+= N-2$ is carried out as in the proof of \rprop {trace2}-Step 3.\qeda

 \medskip  \medskip

 \noindent{\bf Acknowledgements} H. Chen is supported by NSF of China, No: 11726614, 11661045, by the
Jiangxi Provincial Natural Science Foundation, No: 20161ACB20007, and by the Alexander von
Humboldt Foundation.

\end{document}